\documentclass{article}
\usepackage[utf8]{inputenc}
\usepackage{amsmath}
\usepackage{amsfonts}
\usepackage{mathrsfs}
\usepackage{amssymb}
\usepackage{bbm}
\usepackage{bm}
\usepackage{amsthm,mathtools}
\usepackage{dsfont}
\usepackage[normalem]{ulem}
\usepackage{enumerate}
\usepackage{xmpmulti}
\usepackage[english]{babel}
\usepackage[T1]{fontenc}
\usepackage{amsthm, overpic}
\usepackage{verbatim}
\usepackage{tikz}
\usepackage{hyperref}
\newcommand{\bpf}[1][Proof]{{\noindent {\sc #1: }}}
\newcommand{\epf}{{{\hfill $\Box$ \smallskip}}}

\newcommand{\E}{\mathbb{E}}

\newcommand{\N}{\mathbb{N}}
\newcommand{\PP}{\mathbb{P}}

\newcommand{\R}{\mathbb{R}}

\newcommand{\T}{\mathbb{T}}
\newcommand{\Z}{\mathbb{Z}}

\newcommand{\Ebf}{\mathbf{E}}

\newcommand{\Pbf}{\mathbf{P}}

\newcommand{\Pp}{\mathsf{P}}

\newcommand{\Ac}{\mathcal{A}}
\newcommand{\Bc}{\mathcal{B}}
\newcommand{\Dc}{\mathcal{D}}
 
\newcommand{\Fc}{\mathcal{F}}

\newcommand{\Lc}{\mathcal{L}}

\newcommand{\id}{\mathbbm{1}}

\newtheorem{theorem}{Theorem}
\newtheorem{lemma}{Lemma}

\newtheorem{proposition}{Proposition}

\newtheorem{definition}{Definition}

\DeclarePairedDelimiterX{\inner}[2]{\langle}{\rangle}{#1,#2}
\DeclarePairedDelimiterX{\norm}[1]{\|}{\|}{#1}

\title{On global solutions to the semidiscrete stochastic heat equation}
\author{%
  Tobias Hurth \footnote{Department of Mathematics, University of Toronto, Toronto, Canada}%
  \and Konstantin Khanin\footnote{Beijing Institute of Mathematical Sciences and Applications (BIMSA), Beijing, China} 
  \and Beatriz Navarro Lameda \footnote{Department of Mathematics, University College London, London, United Kingdom}%
  }

\date{}

\begin{document} 

\maketitle

\begin{abstract}

We consider the stochastic heat equation on the integer lattice $\Z^d$ in dimension $d \geq 3$ and with small coupling constant. We show uniqueness of global solutions within the class of positive functions that are stationary in time and whose asymptotic growth in space is subexponential. Our proof relies on a factorization formula for the point-to-point partition function in the associated polymer model.  

\end{abstract}

\noindent {\bf Keywords:} Stochastic heat equation, Directed polymers, Weak disorder, Partition function \\
{\bf MSC numbers:} 60H15, 35R60, 37L40, 60K35

\section{Introduction}        \label{sec:intro}

In this paper we address the problem of global solutions to the stochastic heat equation, which is closely related to
the random Hamilton-Jacobi equation. In the introduction below we present a general overview of the area.
We also discuss our main results and their motivation.

The theory of the randomly forced Hamilton-Jacobi equation, and the closely related randomly forced Burgers equation,
has been actively studied in the last 30 years. The initial motivation came from physics, where the main object of
interest was a turbulence-type behavior, so-called ``burgulence" (\cite {Frisch2001, BKh}). Later development was connected
with the Kardar--Parisi--Zhang (KPZ) phenomena and random growth models for interfaces~\cite{Kardar:1986aa}. These topics have been very intensively studied in the last 15 years.

In its most general form the random Hamilton-Jacobi equation can be written in the following way:
\begin{equation}   \label{RFHJ}
\phi_t + H(\nabla \phi) = \nu \Delta \phi + F^\omega,
\end{equation} 
where $\phi$ is a scalar function on space-time $\R^d \times \R$, $H$ is the Hamiltonian, $F^{\omega}$ is a space- and time-dependent random potential, $\nu \geq 0$ is the viscosity, and $\nabla, \Delta$ are the spatial gradient and Laplacian, respectively. The Hamiltonian $H(p)$ is supposed to be convex in the momenta
$p\in \R^d$. It is also assumed that the corresponding Lagrangian $L(v) := \max_p{[p \cdot v - H(p)]}$ grows superlinearly in
$\|v\|$. The random potential $F^{\omega}(x,t)$ is typically stationary in $x$ and $t$. It is also important to assume that its correlations
decay fast in both the space and time variables. It is usually assumed to be white in time; this assumption does not affect the
regularity of the equation. However, a white in space potential creates very serious difficulties on small scales. These difficulties were addressed by M. Hairer through his theory of regularity structures~\cite {Hairer_KPZ} \cite{Hairer_regularity} and by Gubinelli, Imkeller, and Perkowski through their theory of paracontrolled distributions~\cite{gubinelli_imkeller_perkowski} \cite{gubinelli_perkowski}. Since we are interested in large-scale behavior, both in space and time, in order to avoid the regularity issues, we assume that the realizations of the potential $F^{\omega}(x,t)$ are smooth in $x$.

The most important and most studied example of a random Hamilton-Jacobi equation corresponds to the quadratic Hamiltonian $H(p)= \|p\|^2/2$. In this case 
the Lagrangian is also quadratic ($L(v)= \|v\|^2/2$), and the map between Legendre conjugate variables $p$ and $v$ is the identity;
namely, $v(p)=p$. It follows that the velocity field $v=\nabla \phi$ satisfies the random vector Burgers equation:
\begin{equation}   \label{RB}
v_t + (v\cdot\nabla)v = \nu \Delta v + \nabla F^\omega. 
\end{equation}
An unforced version of the Burgers equation was introduced in the 1930s by J. M. Burgers~\cite{zbMATH02512692} as a one-dimensional model for the dynamics of pressure-less gas. He was guided by the observation that the equation has similar invariances, conservation laws, and hydrodynamical nonlinearity as the Navier--Stokes equation, while being easier to study. However, it is the presence of a random forcing term in~\eqref{RB} that describes the physics of turbulence-type behavior~\cite{Frisch2001, BKh}. The randomly forced Burgers equation has been studied extensively in the last 30 years, both in the physics literature~\cite{Krug:1991aa, PhysRevE.51.R2739, PhysRevE.54.4908, MUNIANDY1997352} and in the mathematics literature~\cite{zbMATH01256225, zbMATH00956434, zbMATH00694958, zbMATH00736845, zbMATH00553828, zbMATH00780710, zbMATH00915162, zbMATH01019753, Kifer, zbMATH01322983}. 

The case of a quadratic Hamiltonian is the only case when the nonlinear equation (\ref{RFHJ}) can be linearized.
Using the Hopf--Cole transformation $\phi =-2\nu\log u$~\cite{zbMATH03059136, zbMATH03066118}, we obtain the stochastic heat equation:
\begin{equation}   \label{SH}
u_t  = \nu \Delta u + \frac{1}{2\nu}F^\omega u.
\end{equation} 
The solutions to the stochastic heat equation can be written using the Feynman-Kac formula. One has to introduce directed polymers
which are Wiener processes in the random environment given by the potential $F^\omega$. Then the solution $u(x,t)$ to the Cauchy problem for the
stochastic heat equation can be written in terms of averaging of the initial condition $u(\cdot,0)$ with respect to the probability distribution 
of the polymer endpoint. The directed polymer initiates at point $x$ at time $t$, and  evolves backward in time on the time interval $s\in [0,t]$.

In the case of the quadratic Hamiltonian, all three equations~\eqref{RFHJ},~\eqref{RB},~\eqref{SH} are equivalent, and any asymptotic results about one of them can be
easily reformulated in terms of the the other two.

The compact case of space-periodic potentials $F^{\omega}(x+k,t)=F^{\omega}(x,t), \, k\in \Z^d$, was studied in the series of papers~\cite{zbMATH01256225, EKMS, zbMATH01903736, zbMATH05140623, zbMATH06788737}, resulting in an almost complete understanding of global solutions in this case. These are solutions which can be extended to the whole time interval $t \in \R^1$.
The average velocity $b=\int_{\T^d}{v(x)}dx$, where $\T^d = \R^d / \Z^d$, is the first integral for the Burgers equation~\eqref{RB}.  In the case of the random Hamilton-Jacobi equation and the stochastic heat equation  
the corresponding time-invariant set consists of functions $\phi(x)= b \cdot x +\psi(x)$ and $u(x)=\exp(b \cdot x) f(x)$, respectively, where $\psi(x)$ and $f(x)$ are
periodic: $\psi(x+k)=\psi(x), \, f(x+k)=f(x), \, k\in \Z^d$. One of the main results in the periodic case can be formulated in the following way: 

\smallskip 

Almost surely, for all $b\in \R^d$ and all $\nu\geq 0$, the random Hamilton--Jacobi equation~\eqref{RFHJ} with space-periodic potential admits a unique (up to an additive constant) global solution with average velocity $b$~\cite{zbMATH01256225, zbMATH01903736}. 

\smallskip 

Of course this result implies uniqueness of global solutions to the random Burgers equation and to the stochastic heat equation
(up to a mutiplicative constant). 

In the non-periodic case the situation is much more difficult. In~\cite[Conjecture~1]{zbMATH06867560} the following conjecture was put forward: 

\smallskip 

Almost surely, for all $b \in \R^d$ and all $\nu \geq 0$, the random Hamilton--Jacobi equation~\eqref{RFHJ} with exponentially decaying space-time correlations for $F^{\omega}$ admits a unique (up to an additive constant) time-stationary global solution of the form $\phi_{b, \omega}(x,t) = b \cdot x + \psi_{b,\omega}(x,t)$, where $\psi_{b,\omega}(\cdot,t)$ has sublinear growth for every $t$. 

\smallskip 

At present there are only few mathematical results supporting this conjecture. Uniqueness has been established in spatial dimension $d=1$ for the quadratic Hamiltonian in the viscous ($\nu > 0$) and in the inviscid case ($\nu = 0$), under certain assumptions on the forcing~\cite{R-AS, BCK, BL1, BL2}. Some uniqueness results for global solutions have also been obtained in the so-called weak-disorder case corresponding to spatial dimension $d\geq 3$ and weak potential $F^\omega$, e.g. in~\cite{Kifer}. The notion of weak disorder
reflects the asymptotic behavior of directed polymers which we discussed above in the context of the Feynman-Kac formula. It turns out that  in the limit $t\to \infty$ 
directed polymers can either be localized or diffusive. Localization means that with large probability the endpoint distribution will have bounded variance, while
in the diffusive regime the variance grows linearly in $t$. The case of localization is called the strong-disorder regime, while diffusive behavior corresponds to the
situation of weak disorder. It was shown in~\cite{CarmonaHu, Comets} that strong disorder always takes place in the low dimensional cases $d=1,2$, while in the case $d\geq 3$ there is a transition from weak to strong disorder as the potential $F^\omega(x,t)$ is increasing. The weak-disorder regime was studied starting from the late 1980s.
The diffusive behavior of the directed polymers was rigorously established in the discrete space-time case (\cite{Imbrie1988},\cite{Bolthausen},\cite{Sinai_95}). Kifer~\cite{Kifer} considered the continuous
space-time setting. He proved diffusive behavior of the directed polymers, and deduced uniqueness of global solutions, but in
a rather weak sense. Namely, it was shown that solutions to the Cauchy problem for the stochastic heat equation converge to the unique $L^2$-stationary global solution if the initial condition $u(\cdot,0)$ is itself $L^2$-stationary. This imposes in particular a considerable restriction on the spatial growth of the initial condition. 
Since different global solutions to the random Hamilton--Jacobi equation depend on a parameter $b\in \R^d$ corresponding to a class of functions of the form $\phi(x,t)= b \cdot x +  \psi(x,t)$ where $\psi$ has sublinear growth in $||x||$~\cite{zbMATH06867560, zbMATH01903736},
one really has to establish convergence for all initial conditions $\phi(\cdot, 0)$ with sublinear growth in the case $b=0$. Via the Hopf--Cole transform, the convergence to the global solution to the stochastic heat equation~\eqref{SH} has to be established for all initial conditions $u(\cdot,0)$ with subexponential growth, i.e. $u(\cdot, 0) = \exp w(\cdot)$, where $w$ has sublinear growth. 

Convergence for initial conditions of the above type is the first main result of this paper (Theorem~\ref{thm:pushforward_attraction}). By establishing this type of convergence we are able to prove uniqueness of global solutions to the stochastic heat equation in a very
strong sense: We show that time-stationary global solutions to~\eqref{SH} with subexponential growth in space are unique up to a rescaling at the origin (Theorem~\ref{thm:unique_global_sol}). 
The results proved in (\cite {Imbrie1988},\cite{Bolthausen},\cite{Sinai_95},\cite{Kifer}) are all related to the behavior on the diffusive scale $\|x\|= O(t^{1/2})$. Since we deal with fast growing initial conditions we have
to extend the analysis much further, basically up to the ballistic scale $\|x\|=O(t)$. Such an extension requires a different approach. In particular one has to establish a factorization result
for partition functions corresponding to polymers with endpoints far away from each other. We prove the factorization formula in the separate paper~\cite{HKN}. 
Since we are motivated by the question of uniqueness for the stochastic partial differential equations \eqref{RFHJ}-\eqref{SH}, we work in a continuous-time setting.
At the same time, we work in discrete space because it allows for a more transparent presentation.
As a result we are proving the main theorem
in the case of the semidiscrete stochastic heat equation. We cannot formulate corresponding results for the random Hamilton-Jacobi equation since the Hopf--Cole transformation is not readily available 
in the discrete setting. We certainly believe that our results and methods can be extended to the continuous-space setting, and plan to address this problem in the future.
For simplicity of the presentation we also consider only the case $b=0$. Extension to all $b\in \R^d$ is relatively straightforward and also will be done in the future.

In the model we consider, $F^\omega(x,t)=\dot W^x_t, \, x\in \Z^d$, where $(W^x)_{x \in \Z^d}$ are independent two-sided standard Wiener processes. Independence in $x$ and white behavior in $t$ are
technical conditions which simplify the proofs, but can be weakened by considering weakly dependent potentials. The Gaussianity of the potential  is more important. While we believe that
it is still a technical condition, we use it to prove lower bounds on the partition function. At present all results of this type are based on Talagrand's  approach, and require Gaussian distributions.

The paper has the following structure: In Section~\ref{sec:setting} we describe in detail the semidiscrete stochastic heat equation we consider, as well as the associated directed-polymer model, which gives rise to a global solution (Proposition~\ref{prop:global_stationary_sol}). We also state our main results on convergence of solutions to this global solution in probability (Theorem~\ref{thm:pushforward_attraction}) and on  uniqueness of the global solution (Theorem~\ref{thm:unique_global_sol}). Finally, we present two key ingredients in the proof of convergence to the global solution: a factorization formula for partition functions of polymers with endpoints on a sub-ballistic scale (Section~\ref{ssec:factorization_formula}) and a lower tail estimate of the law for the limiting partition function obtained through Talagrand's method (Section~\ref{ssec:Talagrand}). The proof of Theorem~\ref{thm:pushforward_attraction} is carried out in Sections~\ref{sec:proof_attraction},~\ref{sec:proof_BC_convergence}, and~\ref{sec:p_ratio_uniform_bound}: Section~\ref{sec:proof_attraction} contains the main part. In Section~\ref{sec:proof_BC_convergence} we show that in the limit as $t \to \infty$ the contribution to the solution of the Cauchy problem coming from points $x$ such that $\|x\| > t^{1-\epsilon}$ is negligible for arbitrarily small $\epsilon > 0$. In the proof of Theorem~\ref{thm:pushforward_attraction} this enables us to exploit the above-mentioned factorization formula for partition functions, whose validity is restricted to the sub-ballistic regime. In Section~\ref{sec:p_ratio_uniform_bound}, we prove an estimate on the ratio of transition probabilities for simple symmetric random walk. In Section~\ref{ssec:unique_sd}, we prepare the proof of Theorem~\ref{thm:unique_global_sol} on uniqueness of the global solution by introducing a Markov random dynamical system which describes the evolution in time of solutions to the stochastic heat equation, normalized to equal $1$ at the origin. The proof of Theorem~\ref{thm:unique_global_sol} is then given in Sections~\ref{sec:proof_unique_global_sol},~\ref{sec:proof_unique_sd}, and~\ref{sec:skew_product_proofs}.   

We conclude this introduction with a brief remark on notation. Throughout this article the Euclidean norm and inner product in $\R^d$ are denoted by $\| \cdot \|$ and $\langle \cdot, \cdot \rangle$, respectively. The $1$-norm in $\R^d$ is denoted by $\| \cdot \|_1$. We write $a(t) \lesssim b(t)$ to denote that $a$ is asymptotically less than $b$, i.e., $\limsup_{t\to\infty}{a(t)}/{b(t)} < \infty$.

\vspace{2mm} 

\section*{Acknowledgments}
Part of this paper was written during two-week stays at Mathematisches \\ Forschungszentrum Oberwolfach, in 2018, and Centre International de Rencontres Math\'ematiques (CIRM), in 2019, as part of their respective research in pairs programs. We thank both institutions for their kind hospitality. TH gratefully acknowledges support from the Einstein Foundation Berlin through grant IPF-2021-651(2022-2024) and, until 2020, from the Swiss National Science Foundation through grant 200021 – 175728/1. KK is grateful to acknowledge support from the Natural Sciences and Engineering Research Council of Canada through NSERC Discovery Grant RGPIN-2024-05804.

\vspace{3mm}

\section{Setting and main results}   \label{sec:setting}

\subsection{Setup}

Throughout this paper, we assume that the dimension $d$ of the integer lattice $\Z^d$ is $\geq 3$ and that the coupling constant or inverse temperature $\beta > 0$ is small. In the language of directed-polymer models, this means we consider the regime of weak disorder. Most of the statements in this paper only hold for $\beta$ sufficiently small, i.e., for a given statement there exists $\beta_c > 0$ such that the statement holds for all $\beta \in (0, \beta_c)$. In many instances, we do not know how to express the optimal $\beta_c$ in terms of the given quantities of our model, which is in part related to the fact that we work in continuous time instead of discrete time. 

\subsubsection{The noise}

For $d \geq 3$, let $\Omega$ be the set of functions $\omega : \Z^d \times \R \to \R$ such that for every $x \in \Z^d$, the function $t \mapsto \omega(x,t)$ is continuous and satisfies $\omega(x,0) = 0$. Each $\omega \in \Omega$ represents a realization of the noise in our stochastic model. 
Let $\Fc$ denote the canonical $\sigma$-field on $\Omega$, and let $Q$ be the probability measure on $(\Omega, \Fc)$ under which $(W^x)_{x \in \Z^d}$, defined by $W^x_t(\omega) := \omega(x,t)$, are independent two-sided Wiener processes. To ensure measurability for sets contained in sets of measure zero, we work with the completion of $(\Omega, \Fc, Q)$ (see, e.g., Proposition~3.3.3 in~\cite{dudley_2002}), which we still denote by $(\Omega, \Fc, Q)$. Expectation corresponding to $Q$ will be denoted by $\langle \cdot \rangle$. For any time $s \in \R$, let $\theta_s : \Omega \to \Omega$ be the \emph{Wiener shift} defined by 
$$
\theta_s \left( \omega(x,t) \right) := \omega(x,t+s) - \omega(x,s), 
$$
for all $(x,t) \in \Z^d \times \R$; i.e., every path $\omega(x, \cdot)$ is shifted by $s$ to the left along the time axis and normalized to equal $0$ at time $t=0$. The probability measure $Q$ is invariant with respect to $(\theta_s)_{s \in \R}$, in the sense that for every $s \in \R$ and for every $A \in \Fc$, one has $Q(\theta_s (A)) = Q(A)$. 

\subsubsection{The parabolic Anderson model} 

Let $f: \Z^d \to (0, \infty)$ be a function of subexponential growth and decay; namely, we assume that for a sufficiently small $\epsilon \in (0,1)$, $f$ satisfies
\begin{equation}   \label{eq:asymptotic_exp_growth} 
\lim_{\|x\| \to \infty} \frac{\left \lvert \ln \left(f(x) \right) \right \rvert}{\|x\|^{1-\epsilon}} = 0. 
\end{equation}  
We fix $s \in \R$ and $\beta > 0$, and consider the following Cauchy problem for the semidiscrete stochastic heat equation (sSHE), also known as the \emph{parabolic Anderson model}: 
\begin{equation}
\label{eq:Cauchy_problem_intro}
\tag{PAM}
\left\{
\begin{aligned}
\partial_t u(y,t)  &= \Delta_y u(y,t) + \beta u(y,t) \dot{W}^y_t, \qquad y \in \Z^d, \ t > s, \\
u(y,s) &= f(y), \qquad y \in \Z^d.  	
\end{aligned}
\right.
\end{equation}
Here, $\Delta_y$ is the discrete Laplacian given by 
$$
\Delta_y u(y,t) := \frac{1}{2d} \sum_{z \in \Z^d: \|y-z\|_1 = 1} \left(u(z,t) - u(y,t) \right),
$$
and $\dot{W}_t^y$ is the white noise associated with $W_t^y$. Throughout this paper, we assume that $\beta$ is small. 
We emphasize that most studies of the parabolic Anderson model consider bounded or even localized initial data, whereas the initial data considered in this paper is in a much more general class of functions, namely those with subexponential growth as in~\eqref{eq:asymptotic_exp_growth}.

\subsubsection{Directed polymers in a random potential}  \label{ssec:directed_polymers} 

For any $(x,s) \in \Z^d \times \R$, let $\eta = (\eta_t)_{t \geq s}$ be a continuous-time simple symmetric random walk (SSRW) on $\Z^d$ starting at $\eta_s = x$. The corresponding probability measure is denoted by $\Pbf_{x,s}$ and the corresponding expectation by $\Ebf_{x,s}$. 
We assume that the jumps of $\eta$ occur at random times given by independent exponential clocks; i.e., the times between consecutive jumps form an i.i.d. sequence of exponential random variables with rate $1$. Note that $\eta$ is transient because $d\geq 3$.
If observed over a time interval $[s,t)$, a sample path of $\eta$ (which we shall also denote by $\eta$) is characterized by 
\begin{enumerate}
\item the number $n_{s,t}$ of jumps that occur within the time interval $(s,t)$, \label{jumps}
\item a discrete-time path $\gamma = (\gamma_0, \gamma_1, \ldots, \gamma_{n_{s,t}})$ on $\Z^d$ such that $\gamma_0 = x$ and $\|\gamma_j - \gamma_{j-1}\|_1 = 1$ for $1 \leq j \leq n_{s,t}$, and
\item the jump times $s < s_1 < \ldots < s_{n_{s,t}} <  t$.
\end{enumerate}
It is convenient to introduce the notation $s_0 := s$ and $s_{n_{s,t}+1} := t$, although we do not assume that $s$ and $t$ are jump times. If $s=0$, we will typically write $n_t$ instead of $n_{0,t}$. To a sample path $\eta$ and a realization of the noise $\omega \in \Omega$, we assign the action defined by 
\begin{equation} \label{eq:def_action}
\Ac_s^t(\eta, \omega) := \sum_{j=0}^{n_{s,t}} \left(\omega(\gamma_j, s_{j+1}) - \omega(\gamma_j, s_j) \right).  
\end{equation}
For any time $t > s$ and any site $y \in \Z^d$, denote the probability measure obtained from $\Pbf_{x,s}$ by conditioning on the event $\{\eta_t = y\}$ by $\Pbf_{x,s}^{y,t}$.  
The corresponding expectation is denoted by $\Ebf_{x,s}^{y,t}$. Also set  
$$ 
p_t^y := \Pbf_{0,0}(\eta_t = y) \quad \text{and} \quad q_n^y := \PP(\gamma_n = y \vert \gamma_0 = 0) 
$$ 
for the transition probabilities of a simple symmetric random walk on $\Z^d$ in continuous and discrete time, respectively.

\subsubsection{Partition functions}    \label{ssec:partition_functions} 

For every $\omega \in \Omega$, we define the random normalized \emph{point-to-point partition function} by  
\begin{equation}      \label{eq:point_to_point_part_fun} 
Z_{x,s}^{y,t}(\omega) := e^{-\frac{\beta^2}{2} (t-s)} p_{t-s}^{y-x} \Ebf_{x,s}^{y,t} e^{\beta \Ac_s^t(\cdot, \omega)}
\end{equation} 
if $s < t$, and $Z_{x,s}^{y,s}(\omega) := \id_{x = y}$. 
We also define the random normalized point-to-line partition functions 
\begin{equation}     \label{eq:def_part_fun}  
Z_{x,s}^t(\omega) := \sum_{y \in \Z^d} Z_{x,s}^{y,t}(\omega), \quad \text{and} \quad Z_s^{y,t}(\omega) := \sum_{x \in \Z^d} Z_{x,s}^{y,t}(\omega). 
\end{equation} 
Since $e^{-\frac{\beta^2}{2} (t-s)} \langle e^{\beta \Ac_s^t(\eta, \cdot)} \rangle = 1$ for every $\eta$, these partition functions are normalized in the sense that 
$$
\langle Z_{x,s}^t \rangle = \langle Z_s^{y,t} \rangle = 1.
$$  
Notice that the law of the stochastic process $(Z_{x,s}^{s+t})_{t \geq 0}$ with respect to $Q$ does not depend on $x$ or $s$ because the law for the increments of the Wiener processes $(W^x)_{x \in \Z^d}$ is stationary in space and time, and because the SSRW $\eta$ is homogeneous. 
Besides, $(Z_{x,s}^{s+t})_{t \geq 0}$ and $(Z_{s-t}^{x,s})_{t \geq 0}$ have the same law because of time-reversibility of $\eta$. 

It was shown in~\cite{Carmona_Molchanov} (see also~\cite{zbMATH02099172, zbMATH05081220}) that the unique solution to the Cauchy problem \eqref{eq:Cauchy_problem_intro}, if interpreted as an integral equation in the sense of It\^{o}, is given by  
\begin{equation}\label{eqn:uf_intro}
u^s_f(y,t) 
	:= \sum_{x \in \Z^d} f(x) Z_{x,s}^{y,t}, \qquad t \geq s.  
\end{equation}
Note that $u_f^s(y,s) = f(y)$. In the special case $s = 0$, we usually write $u_f$ instead of $u^0_f$. 
This result can be viewed as a Feynman--Kac formula for the semidiscrete parabolic Anderson problem. In particular, for any $(x,s) \in \Z^d \times \R$, the point-to-point partition function $Z_{x,s}^{y,t}$ solves 
\begin{equation} \label{eq:sSHE} 
\tag{sSHE} 
\partial_t u(y,t)  = \Delta_y u(y,t) + \beta u(y,t) \dot{W}^y_t, \qquad y \in \Z^d, \ t > s. 
\end{equation} 
See also~\cite[Proposition~5]{zbMATH05081220}. 

The following convergence result for the point-to-line partition functions is well known in the discrete-time setting. The continuous-time version stated here is proved in~\cite{HKN}.


\begin{theorem}     \label{thm:limiting_part_fun}
If $\beta$ is sufficiently small, the following statements hold. 
\begin{enumerate} 
\item For all $(x,s), (y,t) \in \Z^d \times \R$, the partition functions $Z_{x,s}^T$ and $Z_S^{y,t}$ converge in $L^2(Q)$ as $T \to \infty$ and $S \to -\infty$, respectively, to the \emph{limiting partition functions} 
$$
	Z_{x,s}^{\infty} := \lim_{T\to \infty} Z_{x,s}^T
\qquad \text{and} \qquad
	Z_{-\infty}^{y,t}  := \lim_{S \to -\infty}  Z_S^{y,t}.
$$
There is $\theta > 0$, independent of $x,y$ and $s,t$, such that 
$$
\lim_{T \to \infty} (T-s)^{\theta} \langle \left(Z_{x,s}^T - Z_{x,s}^{\infty} \right)^2 \rangle = 0 \
\text{and}
\lim_{S \to -\infty} (t-S)^{\theta} \langle \left(Z_S^{y,t} - Z_{-\infty}^{y,t} \right)^2 \rangle = 0. 
$$
\item There is a set $\Omega^{\lim}_+ \in \Fc$ with $Q(\Omega^{\lim}_+) = 1$ and $\theta_t(\Omega^{\lim}_+) = \Omega^{\lim}_+ \ \forall t \in \R$ such that for all $(x,s), (y,t) \in \Z^d \times \R$ and all $\omega \in \Omega^{\lim}_+$, the limiting partition functions $Z_{x,s}^{\infty}(\omega)$ and $Z_{-\infty}^{y,t}(\omega)$ exist as pointwise limits $\lim_{T \to \infty} Z_{x,s}^T(\omega)$ and $\lim_{S \to -\infty} Z_S^{y,t}(\omega)$, respectively, and are positive. 
\end{enumerate} 
\end{theorem} 

Part (1) follows from Theorems~2.1 and~2.2 of~\cite{HKN}, and part (2) follows from Theorem~2.3 of~\cite{HKN}. The limiting partition function $Z_{-\infty}^{y,t}(\omega)$ is of special interest because it defines a global stationary solution to~\eqref{eq:sSHE} (see Proposition~\ref{prop:global_stationary_sol} in Section~\ref{ssec:unique_gss}).

\subsection{Global stationary solutions}      \label{ssec:gss}

In this paper, our main focus is the analysis of solutions to~\eqref{eq:sSHE} that, in the time variable, are both global and stationary.

\begin{definition}
Let $\Omega' \in \Fc$ such that $Q(\Omega') = 1$ and $\theta_t(\Omega') = \Omega'$ for all $t \in \R$. A measurable map $Z: \Z^d \times \R \times \Omega' \to (0, \infty)$ is called a \emph{global stationary solution} to~\eqref{eq:sSHE} if: 
\begin{enumerate}
\item For every $y \in \Z^d$, $s, t \in \R$ with $s < t$, and $\omega \in \Omega'$, 
$$
Z(y,t,\omega) = \sum_{x \in \Z^d} Z(x,s,\omega) Z_{x,s}^{y,t}(\omega); 
$$
\item For every $y \in \Z^d$, $t \in \R$, and $\omega \in \Omega'$, one has $Z(y,t,\omega) = Z(y,0,\theta_t \omega)$; 
\item For every $y \in \Z^d$, the map $\omega \mapsto Z(y,0,\omega)$ is $(\Fc_{-\infty,0}^{'}, \Bc_0^{\infty})$-measurable, where $\Fc_{-\infty,0}^{'}$ is the restriction of $\sigma(W^x_u: u \leq 0, x \in \Z^d)$ to $\Omega'$ and where $\Bc_0^{\infty}$ is the Borel $\sigma$-field on $(0, \infty)$.  
\end{enumerate}
\end{definition}

Notice that positivity is incorporated into our definition of a global stationary solution. A trivial nonpositive map satisfying (1), (2), and (3) of the previous definition is $Z \equiv 0$. 
As mentioned at the end of Section~\ref{ssec:partition_functions}, a particular global stationary solution to~\eqref{eq:sSHE} is given by $Z_{-\infty}^{y,t}(\omega)$. In a strong sense (see Theorem~\ref{thm:unique_global_sol}), this is the only global stationary solution. 

\begin{proposition}    \label{prop:global_stationary_sol} 
For $\beta$ sufficiently small, the following holds: There is an $\Fc$-measurable set $\Omega^{\textup{sol}} \subset \Omega^+$ with $Q(\Omega^{\textup{sol}}) = 1$ such that the function 
\begin{equation}\label{200211143059}
\begin{aligned}
	\Z^d \times \R \times \Omega^{\textup{sol}} 
		&\to (0, \infty)
\\	(y,t,\omega)
		&\mapsto Z_{-\infty}^{y,t}(\omega)
\end{aligned}
\end{equation}
is a global stationary solution to~\eqref{eq:sSHE}. 
\end{proposition}

\bpf By part (2) of Theorem~\ref{thm:limiting_part_fun}, there is $\Omega^{\lim}_+ \in \Fc$ with $Q(\Omega^{\lim}_+) = 1$ such that for every $(y,t) \in \Z^d \times \R$ and $\omega \in \Omega^{\lim}_+$, the limit $Z_{-\infty}^{y,t}(\omega)$ exists in $(0,\infty)$. In addition, by Lemma~8.3 in~\cite{HKN}, there is $\Omega_{\textup{CK}} \in \Fc$ (CK for Chapman--Kolmogorov) with $Q(\Omega_{\textup{CK}}) = 1$ such that for every $(y,t) \in \Z^d \in \R$, $s < t$, and $\omega \in \Omega_{\textup{CK}}$, one has   
$$
Z_{-\infty}^{y,t}(\omega) = \sum_{x \in \Z^d} Z_{x,s}^{y,t}(\omega) Z_{-\infty}^{x,s}(\omega).  
$$
Define 
\begin{align*} 
\Omega^{\textup{sol}} :=& \{\omega \in \Omega: \ \lim_{s \to -\infty} Z_s^{y,t}(\omega) \in (0, \infty) \ \forall (y,t) \in \Z^d \times \R, \\
& Z_{-\infty}^{y,t}(\omega) = \sum_{x \in \Z^d} Z_{-\infty}^{x,s}(\omega) Z_{x,s}^{y,t}(\omega) \ \forall (y,t) \in \Z^d \times \R, \ \forall s < t\}. 
\end{align*} 
Then $\Omega^{\lim}_+ \cap \Omega_{\textup{CK}} \subset \Omega^{\textup{sol}}$. As $Q(\Omega^{\lim}_+ \cap \Omega_{\textup{CK}}) = 1$ and as $(\Omega, \Fc, Q)$ is complete, the set $\Omega^{\textup{sol}}$ is measurable and has measure $1$. It is straightforward to check that $\theta_t(\Omega^{\textup{sol}}) = \Omega^{\textup{sol}}$ for every $t \in \R$ and that $Z_{-\infty}^{y,t}(\omega) = Z_{-\infty}^{y,0}(\theta_t \omega)$ for every $(y,t) \in \Z^d \times \R$ and $\omega \in \Omega^{\textup{sol}}$. The measurability of $\Z^d \times \R \times \Omega^{\textup{sol}} \ni (y,t,\omega) \mapsto Z_{-\infty}^{y,t}(\omega)$ is not hard to verify either. Finally, it is clear that for every $y \in \Z^d$, the function $\omega \mapsto Z_{-\infty}^{y,0}(\omega)$ is $(\Fc_{-\infty,0}^{\textup{sol}}, \Bc_0^{\infty})$-measurable, with $\Fc_{-\infty,0}^{\textup{sol}}$ the restriction of $\sigma(W^x_u: u \leq 0, x \in \Z^d)$ to $\Omega^{\textup{sol}}$. Then the definition of $\Omega^{\textup{sol}}$ implies that $\Z^d \times \R \times \Omega^{\textup{sol}} \ni (y,t,\omega) \mapsto Z_{-\infty}^{y,t}(\omega)$ is a global stationary solution to~\eqref{eq:sSHE}.
\epf 

\bigskip 

\subsection{Attraction to the global solution}

The first main result of this paper is Theorem~\ref{thm:pushforward_attraction} below, which states that the particular global stationary solution $Z_{-\infty}^{y,t}$ from~\eqref{200211143059} attracts solutions to the Cauchy problem \eqref{eq:Cauchy_problem_intro} with any subexponentially growing initial data $f$. As we shall see in Section~\ref{ssec:unique_gss}, this implies the uniqueness of the global stationary solution in a strong sense.  
For $c > 0$ and $\epsilon \in (0,1]$, let $\Lc_{c,\epsilon}$ be the set of functions $f: \Z^d \to (0,\infty)$ such that 
\begin{equation}   \label{eq:sublinear_log} 
\left \lvert \ln \left(f(x) \right)  \right \rvert \leq c \|x\|^{1-\epsilon}, \qquad \forall x \in \Z^d. 
\end{equation}
Note that this condition implies $f(0) = 1$ and is equivalent to 
\begin{equation}    \label{eq:f_bounds} 
e^{-c \|x\|^{1-\epsilon}} \leq f(x) \leq e^{c \|x\|^{1-\epsilon}}, \qquad \forall x \in \Z^d. 
\end{equation} 

\begin{theorem}     \label{thm:pushforward_attraction}
For $\beta$ sufficiently small, the following holds: For every $y \in \Z^d$ and for every $c > 0$, $\epsilon \in (0,1)$, and $f \in \Lc_{c, \epsilon}$, one has   %
\begin{equation} \label{eq:main_thm}
\left \lvert \frac{u_f(y,t)}{u_f(0,t)} - \frac{Z_{-\infty}^{y,t}}{Z_{-\infty}^{0,t}} \right \rvert \xrightarrow[t \to \infty] {} 0  \text{ in probability. } 
\end{equation}
\end{theorem}

The proof of Theorem~\ref{thm:pushforward_attraction} is given in Section~\ref{sec:proof_attraction}. 

\subsection{Uniqueness of the global stationary solution}   \label{ssec:unique_gss}

The second main result of this paper, and a consequence of Theorem~\ref{thm:pushforward_attraction}, is that, up to a random multiplicative constant, the global stationary solution $Z_{-\infty}^{y,t}(\omega)$ is unique within the class of subexponentially growing functions, implying in particular that the solutions to \eqref{eq:Cauchy_problem_intro} are asymptotically independent of the initial data.

\begin{theorem}  \label{thm:unique_global_sol}
If $\beta$ is sufficiently small, the following holds: Let $\Omega' \in \Fc$ such that $Q(\Omega') = 1$ and $\theta_t(\Omega') = \Omega'$ for all $t \in \R$. Let $Z: \Z^d \times \R \times \Omega' \to (0, \infty)$ be a global stationary solution to~\eqref{eq:sSHE} such that for every $\omega \in \Omega'$ there is $\epsilon \in (0,1)$ with 
$$
\lim_{\|y\| \to \infty} \frac{\lvert \ln(Z(y,0,\omega)) \rvert}{\|y\|^{1-\epsilon}} = 0. 
$$
Then, $Q$-almost surely, for all $(y,t) \in \Z^d \times \R$, 
$$
Z(y,t,\omega) = Z_{-\infty}^{y,t}(\omega) \cdot \frac{Z(0,0,\omega)}{Z_{-\infty}^{0,0}(\omega)}. 
$$
\end{theorem} 

Theorem~\ref{thm:unique_global_sol} is proved in Section~\ref{sec:proof_unique_global_sol}. 
Let us make a few remarks regarding our assumptions.
\begin{enumerate}
\item As it was already mentioned before, we cannot formulate a corresponding result for the random Hamilton-Jacobi equation \eqref{RFHJ} or the random Burgers equation \eqref{RB} because the Hopf--Cole transformation is not readily available in discrete space.
However, we believe that the results and methods of this paper can be extended to the continuous-space setting.
\item We also note that for simplicity of presentation we have only considered the case of global solutions with subexponential growth, which corresponds to the case $b = 0$ for the average velocity introduced in Section~\ref{sec:intro}. We expect the extension of our results to the general case $b \in \R^d$ to be relatively straightforward. In continuous space, this extension would be a simple consequence of shear-invariance, whereas shearing in discrete space will require minor technical adjustments. 
\item The random potential $F^{\omega}$ in our model is the white noise associated with space-independent two-sided standard Wiener processes.
The independence in $x$ and the white behavior in $t$ are purely technical conditions that simplify our proofs; they can be relaxed by considering weakly dependent potentials. The Gaussianity of the potential is more important because, in Theorem~\ref{thm:continuous_Talagrand}, the lower tail estimate of the probability distribution for the point-to-line partition function is based on Talagrand's approach, which requires Gaussian distributions.
Nevertheless, we believe that it is still a technical condition. 
We should mention that lower tail estimates in the general case were recently obtained by St. Junk and H. Lacoin 
(\cite{St_Junk_Hubert_Lacoin}).
\item
Although our results are formulated for $\beta$ small enough, we believe that they can be extended to the whole of the $L^2$ regime. In fact, recent results by St. Junk give hope to extend the main results even further.
Namely, to all $\beta$
corresponding to the weak disorder behavior (\cite{St_Junk}).

\end{enumerate}

\subsection{Factorization formula}    \label{ssec:factorization_formula}

The proof of Theorem~\ref{thm:pushforward_attraction} relies on a factorization formula for the point-to-point partition function $Z_{x,s}^{y,t}$, obtained in~\cite{HKN}. In the high dimension, high temperature regime (i.e., $d \geq 3$ and small $\beta$), Sinai~\cite{Sinai_95}, Kifer~\cite{Kifer}, and Vargas~\cite{Vargas} established a factorization formula for $Z_{x,s}^{y,t}$ in different polymer models, which can be viewed as a local limit theorem for directed polymers in a random environment. Their results pertain to the behavior on the diffusive scale $\|y- x \| = O((t-s)^{1/2})$. In~\cite{Kifer}, Kifer used this kind of factorization formula to show that solutions to the Cauchy problem for the stochastic heat equation converge to the unique space-time stationary global solution if the initial condition is $L^2$ stationary. However, since we are dealing with fast-growing initial conditions, it becomes necessary to extend the analysis far beyond the diffusive scale, essentially up to the ballistic scale $\| y - x \| = O(t-s)$, i.e., to polymers whose starting and endpoint are far away from each other. The following theorem was proved in~\cite[Theorem~2.4]{HKN}. 
\begin{theorem}     \label{thm:factorization}
For $\beta$ sufficiently small, the following holds: For every $\sigma \in (0,1)$ there exists $\theta = \theta(\sigma) > 0$ such that for all $x,y\in\Z^d$ and $s < t$ with $\| x - y \| < (t - s)^{\sigma}$, the partition function $Z_{x,s}^{y,t}$ has the representation 
\begin{equation}     \label{eq:factor_formula}
Z_{x,s}^{y,t} = p_{t-s}^{y-x} \left(Z_{x,s}^{\infty} Z_{-\infty}^{y,t} + \delta_{x,s}^{y,t} \right),
\end{equation} 
where the error term $\delta_{x,s}^{y,t}$ defined by the formula above satisfies 
\begin{equation}   \label{eq:convergence_error}
\lim_{(t-s) \to \infty} (t-s)^{\theta} \sup_{x, y \in \Z^d: \|x-y\| < (t-s)^{\sigma}} \langle \lvert \delta_{x,s}^{y,t} \rvert \rangle = 0.
\end{equation}
\end{theorem}

Intuitively, this formula says that, even if starting point $x$ and endpoint $y$ are far away from each other, the partition function $Z_{x,s}^{y,t}$ only ``feels" the environment at times close to $s$, when it stays near $x$, and at times close to $t$, when it stays near $y$.  

\subsection{A lower tail estimate of the probability distribution for the partition function: Talagrand's method}    \label{ssec:Talagrand}

Another key element of the proof of Theorem~\ref{thm:pushforward_attraction} is a lower tail estimate of the probability distribution for the partition function $Z_0^{y,t}$, which was also obtained in~\cite{HKN}.
In the case of discrete space-time, such estimates have been derived by P. Carmona and Hu \cite{CarmonaHu} using concentration of measure arguments for discrete directed polymers in Gaussian environments that originated in Talagrand's work on spin glasses \cite{zbMATH01129463,Talagrand_11}.
As we work in semidiscrete space-time, we have to adapt the argument by Carmona and Hu to our setting. The theorem below is the continuous-time version of \cite[Theorem~1.5]{CarmonaHu} (see also \cite[Theorem 1(a)]{zbMATH06333766}). It is proved in~\cite[Theorem~2.5]{HKN}. 

\begin{theorem}       \label{thm:continuous_Talagrand}
For $\beta$ sufficiently small, there exists a constant $c > 0$ such that 
\begin{equation*}
Q \left(Z_0^{y,t} < e^{-u} \right) < c e^{-u^2 / c}, \qquad \forall t, u > 0. 
\end{equation*}
\end{theorem}

Theorem~\ref{thm:continuous_Talagrand} implies that the limiting partition function $Z_{-\infty}^{0,0}$ admits all negative moments.

\vspace{3mm}

\section{Proof of Theorem~\ref{thm:pushforward_attraction}}    \label{sec:proof_attraction} 

In this section we prove Theorem~\ref{thm:pushforward_attraction} on convergence in probability of ratios for solutions $u_f$ to the Cauchy problem in~\eqref{eq:Cauchy_problem_intro} to ratios for $Z_{-\infty}^{y,t}$. We start with a brief sketch of the proof strategy.

For given $c > 0$ and $\epsilon \in (0,1)$, the first step is to fix $\sigma \in (1/(1+\epsilon), 1)$, and, for $f \in \Lc_{c, \epsilon}$, to rewrite the solution $u_f$ as the sum of two terms
\begin{equation}     \label{eq:u_sum} 
u_f(y,t) = \sum_{\|x\| \leq t^{\sigma}} f(x) Z_{x,0}^{y,t} + \sum_{\|x\| > t^{\sigma}} f(x) Z_{x,0}^{y,t}
\end{equation} 
corresponding to terms inside and outside of the ball of radius $t^\sigma$, respectively.
The contribution of the second term to $u_f(y,t)$ is negligible: 

One has 
\begin{align*}
\biggl\langle \sum_{\|x\| > t^{\sigma}} f(x) Z_{x,0}^{y,t} \biggr\rangle \leq& \biggl\langle \sum_{\|x\| > t^{\sigma}} e^{c \|x\|^{1-\epsilon}} Z_{x,0}^{y,t} \biggr\rangle \\
=& \sum_{\|x\| > t^{\sigma}} e^{c \|x\|^{1-\epsilon}} p_t^{y-x} \lesssim \sum_{\|x\| > t^{\sigma}} e^{c \|x\|^{1-\epsilon} - \kappa \frac{\|x\|^2}{t}},
\end{align*}  
for some $\kappa > 0$. The expression on the right-hand side is negligible because $\sigma > \frac{1}{1+\epsilon}$. 
The dominant contribution to $u_f(y,t)$ then comes from the first term on the right-hand side of~\eqref{eq:u_sum}. To deal with it, we apply the factorization formula for the partition function $Z_{x,s}^{y,t}$ from Theorem~\ref{thm:factorization}. Using the representation in \eqref{eq:factor_formula}, we have
\begin{align} \label{eq: sum_inside}
\sum_{\|x\| \leq t^{\sigma}} f(x) Z_{x,0}^{y,t} 
&= \sum_{\|x\| \leq t^{\sigma}} f(x) \ p_t^{y-x} \left(Z_{x,0}^{\infty} Z_{-\infty}^{y,t} + \delta_{x,0}^{y,t} \right)\\ \notag
&=\sum_{\|x\| \leq t^{\sigma}} f(x) \ p_t^{y-x} Z_{x,0}^{\infty} Z_{-\infty}^{y,t} + \sum_{\|x\| \leq t^{\sigma}} f(x) \ p_t^{y-x} \delta_{x,0}^{y,t}.
\end{align}
Our aim is to show that the first term in the second line of~\eqref{eq: sum_inside} dominates the sum. To prove this fact, we use the smallness of $\delta_{x,0}^{y,t}$ implied by~\eqref{eq:convergence_error}. However, we must also make sure that $Z_{-\infty}^{y,t}$ does not become too small for typical realizations of the noise. 
For this the main ingredient is Theorem~\ref{thm:continuous_Talagrand}, which implies that for $\theta > 0$, with high probability, 
$$ 
Z_{-\infty}^{y,t} Z_{x,0}^{\infty} \geq t^{-\theta}, \quad \forall x: \|x\| \leq t^{\sigma}. 
$$
This allows us to conclude that the second term in~\eqref{eq: sum_inside} is indeed negligible compared to the first, and therefore
$$
\frac{u_f(y,t)}{u_f(0,t)} \approx \frac{Z_{-\infty}^{y,t} \sum_{\|x\| \leq t^{\sigma}} f(x) \ p_t^{y - x} Z_{x,0}^{\infty}}{Z_{-\infty}^{0,t} \sum_{\|x\| \leq t^{\sigma}} f(x) \ p_t^{x} Z_{x,0}^{\infty}}.
$$ 
Finally, for large $t$, the ratio ${p_t^{y -x}}/{p_t^{x}}$ is close to $1$, uniformly in $\|x\| \leq t^{\sigma}$ (see Lemma~\ref{lm:p_ratio_uniform_bound} below), so the right-hand side is roughly ${Z_{-\infty}^{y,t}}/{Z_{-\infty}^{0,t}}$, which is what we want to show. 
\vspace{3mm}

We now give a formal proof of Theorem~\ref{thm:pushforward_attraction}.
Fix $c > 0$, $\epsilon \in (0,1)$, $f \in \Lc_{c,\epsilon}$, and $\sigma \in (1/(1+\epsilon),1)$. Note that $\sigma > 1/2$ because $\epsilon < 1$. For any $y \in \Z^d$ and $t > 0$, we use the factorizarion formula from Theorem~\ref{thm:factorization} to write $u_f(y,t)$ as the sum of three terms, which we denote by $\Sigma^{(1)}(y,t), \Sigma^{(2)}(y,t)$, and $\Sigma^{(3)}(y,t)$:
\begin{align*} 
u_f(y,t) &=  \sum_{\|x\| \leq t^{\sigma}} f(x)  p_t^{y-x} Z_{x,0}^{\infty} Z_{-\infty}^{y,t} + \sum_{\|x\| \leq t^{\sigma}} f(x) p_t^{y-x} \delta_{x,0}^{y,t}
 + \sum_{\|x\| > t^{\sigma}} f(x) Z_{x,0}^{y,t} \\
& =: \Sigma^{(1)}(y,t) + \Sigma^{(2)}(y,t) + \Sigma^{(3)}(y,t). 
\end{align*}

\noindent Then we factorize $u_f$ as follows:
\begin{equation*}
u_f = \Sigma^{(1)} \left(1 + \frac{\Sigma^{(2)}}{\Sigma^{(1)}}\right)\left(1 + \frac{\Sigma^{(3)}}{\Sigma^{(1)} + \Sigma^{(2)}}\right).
\end{equation*}
By Theorem~\ref{thm:limiting_part_fun}, there is no danger of dividing by $0$, except on a set of measure zero. 
Let 
$$
B(y,t) := \frac{\Sigma^{(2)}(y,t)}{\Sigma^{(1)}(y,t)}, \quad C(y,t) := \frac{\Sigma^{(3)}(y,t)}{\Sigma^{(1)}(y,t) + \Sigma^{(2)}(y,t)}, \quad D(y,t) := \frac{\Sigma^{(1)}(y,t)}{Z^{y,t}_{-\infty}}. 
$$
Then 
$$
\frac{u_f(y,t)}{u_f(0,t)} - \frac{Z^{y,t}_{-\infty}}{Z^{0,t}_{-\infty}} = \frac{Z^{y,t}_{-\infty}}{Z^{0,t}_{-\infty}} \biggl(\frac{(1+B(y,t)) \ (1+C(y,t)) \ D(y,t)}{(1+B(0,t)) \ (1+C(0,t)) \ D(0,t)} - 1 \biggr). 
$$
Since the law of $Z^{y,t}_{-\infty}/Z^{0,t}_{-\infty}$ is a probability measure on $(0, \infty)$ which does not depend on $t$, it is enough to prove that 
\begin{equation*} 
\biggl \lvert \frac{(1 + B(y,t)) \ (1+C(y,t)) \ D(y,t)}{(1+B(0,t)) \ (1+C(0,t)) \ D(0,t)} - 1 \biggr \rvert \xrightarrow[t \to \infty] {} 0  \text{ in probability. } 
\end{equation*}
To deal with the ratio $D(y,t)/D(0,t)$, we use the following result proved in Section~\ref{sec:p_ratio_uniform_bound}.  

\begin{lemma}     \label{lm:p_ratio_uniform_bound}
Let $\sigma \in (0,1)$ and $y_1, y_2 \in \Z^d$. Then, 
$$ 
\lim_{t \to \infty} \inf_{\|x\| \leq t^{\sigma}} \frac{p_t^{y_1-x}}{p_t^{y_2-x}} = \lim_{t \to \infty} \sup_{\|x\| \leq t^{\sigma}} \frac{p_t^{y_1-x}}{p_t^{y_2-x}} = 1. 
$$ 
\end{lemma} 

By Lemma~\ref{lm:p_ratio_uniform_bound}, for every $\zeta > 0$ there exists $T > 0$ such that $1 - \zeta < p^{y-x}_t/p^x_t < 1 + \zeta$ for all $t \geq T$ and $\|x\| \leq t^{\sigma}$. Hence, 
$$
1 - \zeta < \frac{D(y,t)}{D(0,t)} < 1 + \zeta 
$$
for all $t \geq T$. Notice that $T$ only depends on $\zeta$ and $y$ and is in particular independent of the noise. This shows that 
$$
\lim_{t \to \infty} \biggl \lvert \frac{D(y,t)}{D(0,t)} - 1 \biggr \rvert = 0. 
$$
Theorem~\ref{thm:pushforward_attraction} will then follow once we prove that $B(y,t)$ and $C(y,t)$ converge to $0$ in probability.

\begin{lemma}      \label{lm:BC_convergence} 
Let $\beta$ be so small that the conclusions of Theorems~\ref{thm:factorization} and~\ref{thm:continuous_Talagrand} hold. Then for $y$ fixed we have the following convergence statements:
\begin{equation} \label{eq: convergence_B}
\lvert B(y,t) \rvert \xrightarrow[t \to \infty]{} 0 \; \text{ in probability},
\end{equation}
\begin{equation} \label{eq: convergence_C}
C(y,t)  \xrightarrow[t \to \infty]{} 0  \; \text{ in probability}.
\end{equation}
\end{lemma}

\noindent Section~\ref{sec:proof_BC_convergence} is devoted to the proof of this lemma. 

\section{Proof of Lemma \ref{lm:BC_convergence} -- Small Contributions} \label{sec:proof_BC_convergence}

\subsection{Proof of ~\eqref{eq: convergence_B}}

By Theorem~\ref{thm:factorization}, there is $\theta > 0$ such that 
$$
\lim_{t \to \infty} t^{\theta} \sup_{x: \|x\| \leq t^{\sigma}} \langle \lvert \delta_{x,0}^{y,t} \rvert \rangle = 0. 
$$
For $t > 0$, define the event 
$$ 
A(t) := \bigcup_{x: \|x\| \leq t^{\sigma}} \{\omega \in \Omega^{\lim}_+: \ Z_{-\infty}^{y,t}(\omega) Z_{x,0}^{\infty}(\omega) < t^{-\theta}\},  
$$ 
where $\Omega^{\lim}_+$ was introduced in Theorem~\ref{thm:limiting_part_fun}. 
Theorem~\ref{thm:continuous_Talagrand} implies that there is $c > 0$ such that 
$$
Q(Z_{0,0}^{\infty} < e^{-u}) \leq c e^{-u^2/c}, \quad \forall u > 0.
$$
Thus, 
$$ 
Q(A(t)) \leq \sum_{x: \|x\| \leq t^{\sigma}} Q(Z_{-\infty}^{y,t} Z_{x,0}^{\infty} < t^{-\theta}) \lesssim t^{d \sigma} Q(Z_{0,0}^{\infty} < t^{-\frac{\theta}{2}}) \lesssim t^{d \sigma - \frac{\theta^2}{4c} \ln(t)}, 
$$  
which tends to $0$ as $t \to \infty$. Fix $\rho > 0$ and let $\tau > 0$ be so large that for every $t \geq \tau$, 
$$ 
Q(A(t)) < \rho
\quad \text{ and } \quad
t^{\theta} \sup_{\|x\| \leq t^{\sigma}} \langle \lvert \delta_{x,0}^{y,t} \rvert \rangle < \rho^2. 
$$ 
If $A(t)^c$ denotes the complement of $A(t)$ in $\Omega^{\lim}_+$, we have for $t \geq \tau$ the estimate
\begin{align*}
Q(\lvert B(y,t) \rvert > \rho) =& Q(\lvert B(y,t) \rvert > \rho, A(t)) + Q(\lvert B(y,t) \rvert > \rho, A(t)^c) \\
<& \frac{\rho}{2} + Q \biggl(A(t)^c, \rho < \sum_{\|x\| \leq t^{\sigma}} \frac{f(x) p_t^{y-x} \lvert \delta_{x,0}^{y,t} \rvert}{\sum_{\|z\| \leq t^{\sigma}} f(z) p_t^{y-z} Z_{z,0}^{\infty} Z_{-\infty}^{y,t}} \biggr) \\
\leq& \frac{\rho}{2} + Q \biggl(\rho < \sum_{\|x\| \leq t^{\sigma}} \frac{f(x) p_t^{y-x}}{\sum_{\|z\| \leq t^{\sigma}} f(z) p_t^{y-z}} t^{\theta} \lvert \delta^{y,t}_{x,0} \rvert \biggr) \\
\leq& \frac{\rho}{2} + \rho^{-1} \sum_{\|x\| \leq t^{\sigma}} \frac{f(x) p_t^{y-x}}{\sum_{\|z\| \leq t^{\sigma}} f(z) p_t^{y-z}} t^{\theta} \langle \lvert \delta_{x,0}^{y,t} \rvert \rangle \\
\leq& \frac{\rho}{2} + \rho^{-1} t^{\theta} \sup_{\|x\| \leq t^{\sigma}} \langle \lvert \delta_{x,0}^{y,t} \rvert \rangle < \rho, 
\end{align*} 
where we used Markov's inequality.

\vspace{-5mm}
\hfill $\qed$



\subsection{Proof of~\eqref{eq: convergence_C}}

\begin{lemma}\label{lm:sum_of_hZ}
One has 
$$
\lim_{t\to\infty} \left\langle \Biggl \lvert \sum_{\|x\| > t^{\sigma}} e^{c \|x\|^{1-\epsilon}} Z_{x,0}^{y,t} \Biggr\lvert \right \rangle = 0.
$$
\end{lemma}

\bpf It is easy to see from the definition of the partition function $Z_{x,0}^{y,t}$ that $\left\langle Z_{x,0}^{y,t} \right\rangle = p_t^{y-x}$. Therefore, 
\begin{equation}    \label{eq:sum_of_hZ_1} 
 \left\langle \biggl \lvert \sum_{\|x\| > t^{\sigma}} e^{c \|x\|^{1-\epsilon}} Z_{x,0}^{y,t} \biggr \rvert \right \rangle = \sum_{\|x\| > t^{\sigma}} e^{c \|x\|^{1-\epsilon}} p_t^{y-x}. 
\end{equation} 
Recall that the continuous-time transition probability $p_{t}^{y-x}$ satisfies
\begin{equation}\label{eq:cts-time_tp}
p_{t}^{y-x} = \sum_{n=0}^\infty e^{-t} \frac{t^n}{n!} q_n^{y-x}, 
\end{equation}
where $q_j^z := \PP(\gamma_j = z \vert \gamma_0 = 0)$ is the transition probability for a discrete-time simple symmetric random walk $(\gamma_j)_{j \in \N_0}$ on $\Z^d$. If $t$ is sufficiently large, we obtain the estimate
\begin{equation}   \label{eq:y_x_estimate_large_t} 
\|y- x \|_1^2 \geq \| y - x\|^2 \geq \frac{1}{2} \| x \|^2 > \frac{1}{2} t^{2\sigma} \quad \text{for all } x \text{ such that } \| x \| > t^\sigma, 
\end{equation} 
so $q^{y-x}_n = 0$ for all $n \leq 2^{-1/2} t^{\sigma}$ and $x \in \Z^d$ such that $\|x\| > t^{\sigma}$. Furthermore, for any fixed $n > 2^{-1/2} t^{\sigma}$, we have $q^{y-x}_n = 0$ if $\|x\| \geq 2n$. Using these two observations together with~\eqref{eq:cts-time_tp} we can rewrite the righthand side of~\eqref{eq:sum_of_hZ_1} as 
\begin{equation}   \label{eq:term_expansion} 
\sum_{n > 2^{-1/2} t^{\sigma}} e^{-t} \frac{t^n}{n!} \sum_{t^{\sigma} < \|x\| < 2n} e^{c \|x\|^{1-\epsilon}} q_n^{y-x}. 
\end{equation} 
Using the estimate 
$$
q_j^z \leq \PP \left(\max_{0 \leq i \leq j} \| \gamma_i \| \geq \|z\| \right), 
$$
and Proposition~2.1.2 in~\cite{Lawler_Limic}, one has  
\begin{equation}   \label{eq:large_deviations} 
q_j^z \leq \hat c e^{-\kappa \|z\|^2 / j}, \quad j \in \N, \ z \in \Z^d, 
\end{equation} 
for some constants $\hat c, \kappa > 0$. 
On account of~\eqref{eq:large_deviations} and~\eqref{eq:y_x_estimate_large_t}, the expression in~\eqref{eq:term_expansion}, for large $t$, is bounded above by a constant times 
$$
\sum_{n > 2^{-1/2} t^{\sigma}} e^{-t} \frac{t^n}{n!} \sum_{t^{\sigma} < \|x\| < 2n} e^{c \|x\|^{1-\epsilon} - \kappa \frac{\|x\|^2}{2n}}. 
$$
Fix $\xi \in (1, \sigma (1+\epsilon))$ and split the sum above into 
\begin{equation}   \label{eq:Dickensian_sums} 
\sum_{2^{-1/2} t^{\sigma} < n \leq t^{\xi}} e^{-t} \frac{t^n}{n!} Y_n(t) + \sum_{n > t^{\xi}} e^{-t} \frac{t^n}{n!} Y_n(t), 
\end{equation} 
where 
$$
Y_n(t) := \sum_{t^{\sigma} < \|x\| < 2n} e^{c \|x\|^{1-\epsilon} - \kappa \frac{\|x\|^2}{2n}}. 
$$
For $n \leq t^{\xi}$ and $x \in \Z^d$ such that $\|x\| > t^{\sigma}$, we have 
$$
c \|x\|^{1-\epsilon} - \kappa \frac{\|x\|^2}{2n} \leq \|x\|^2 \left(c t^{-\sigma (1+\epsilon)} - \frac{\kappa}{2} t^{-\xi} \right) \leq -\|x\|^2 \frac{\kappa}{4} t^{-\xi}
$$
provided that $t$ is sufficiently large.  
This yields for $n \leq t^{\xi}$ and $t$ large enough 
$$
Y_n(t) \leq \sum_{t^{\sigma} < \|x\| < 2n} e^{-\frac{\kappa}{4}\frac{\|x\|^2}{t^{\xi}}}.
$$
Since 
$$
\bigl \lvert \{x \in \Z^d: \|x\| \leq r\} \bigr\lvert = O(r^d),  
$$
we have for $n \leq t^{\xi}$ the estimate 
$$ 
\bigl \lvert \left\{x \in \Z^d: t^{\sigma} < \|x\| < 2n \right\} \bigr\lvert \lesssim t^{\xi d}. 
$$
As a result, 
\begin{equation}   \label{eq:Y_series_estimate}
Y_n(t) \lesssim t^{\xi d} e^{-\frac{\kappa}{4} t^{2 \sigma - \xi}} \xrightarrow[t \to \infty] {} 0. 
\end{equation}
Since the upper bound in~\eqref{eq:Y_series_estimate} does not depend on $n$, we also have 
$$
\lim_{t\to\infty} \sum_{2^{-1/2} t^{\sigma} < n \leq t^{\xi}} e^{-t} \frac{t^n}{n!} Y_n(t) =0.
$$
To deal with the second sum in~\eqref{eq:Dickensian_sums}, fix $\zeta  > 1$ and notice that 
$$
Y_n(t) \leq \sum_{t^{\sigma} < \|x\| < 2n} e^{c \|x\|^{1-\epsilon}} \lesssim n^d e^{c (2n)^{1-\epsilon}} \lesssim \zeta^n. 
$$
The second sum in~\eqref{eq:Dickensian_sums} is therefore bounded from above by a constant times 
\begin{equation}    \label{eq:Y_1_estimate_1}
\sum_{n > t^{\xi}} e^{-t} \frac{(\zeta t)^n}{n!}.  
\end{equation}
Using the tail estimate 
$$
\sum_{n=k}^{\infty} \frac{s^n}{n!} \leq \frac{s^k}{k!} \sum_{n=k}^{\infty} \left(\frac{s}{k} \right)^{n-k} = \frac{s^k}{k!} \frac{1}{1-\frac{s}{k}}, \quad k > s 
$$
and Stirling's formula, we see that~\eqref{eq:Y_1_estimate_1} is bounded from above by a constant times 
$$
t^{-\frac{\xi}{2}} e^{-t} \left(\frac{e \zeta t}{\lfloor t^{\xi} \rfloor} \right)^{\lfloor t^{\xi} \rfloor} \xrightarrow[t \to \infty] {} 0.
$$  
This completes the proof. 
\epf 
\bigskip

\begin{lemma}       \label{lm:relevant_denominator}
Let $\beta$ be so small that the conclusion of Theorem~\ref{thm:continuous_Talagrand} holds. Then, for every $\delta > 0$ there exist $u>0$ and $T_1>0$ such that for all $t \geq T_1$, 
$$
Q \Biggl( \sum_{\|x\|\leq t^{\sigma}} Z_{x,0}^{y,t} < e^{-u} \Biggr)< \delta.
$$ 
\end{lemma}

\bpf For $u>2\ln 2$, we have  
\begin{align}    \label{eq:pf_relevant_denominator} 
Q\Biggl( \sum_{\|x\|\leq t^{\sigma}} Z_{x,0}^{y,t} < e^{-u} \Biggr)
	&= Q \Biggl(Z_0^{y,t} - \sum_{\|x\|>t^{\sigma}} Z_{x,0}^{y,t} < e^{-u} \Biggr)  \\
	&\leq Q\left(Z_0^{y,t} < 2 e^{-u}\right) + Q\Biggl( \sum_{\|x\|>t^{\sigma}} Z_{x,0}^{y,t} > e^{-u} \Biggr) \notag \\
&\leq Q \left(Z_0^{y,t} < e^{-\frac{u}{2}} \right) + Q \Biggl(\sum_{\|x\| > t^{\sigma}} Z_{x,0}^{y,t} > e^{-u} \Biggr).  \notag
\end{align}
Theorem~\ref{thm:continuous_Talagrand} and Markov's inequality imply that the third line of~\eqref{eq:pf_relevant_denominator} is less than 
$$ 
ce^{- \frac{u^2}{4c}} + e^u \left \langle \sum_{\|x\| > t^{\sigma}} Z_{x,0}^{y,t} \right\rangle = ce^{-\frac{u^2}{4c}} + e^u \sum_{\|x\| > t^{\sigma}} p_t^{y-x}
$$ 
for some $c > 0$ that does not depend on $u$ or $t$. Fix $\delta > 0$ and let $u$ be so large that $ce^{-\frac{u^2}{4c}} < \frac{\delta}{2}$. Since $\sum_{\|x\| > t^{\sigma}} p_t^{y-x} \longrightarrow 0$ as $t \to \infty$, there exists $T_1>0$ such that for all $t \geq T_1$, we have $e^u \sum_{\|x\| > t^{\sigma}} p_t^{y-x} <  \frac{\delta}{2}$.

\epf

\bigskip

\bpf[Proof of ~\eqref{eq: convergence_C}] Because of~\eqref{eq:f_bounds} we have   
\begin{equation*}
C(y,t)  = \frac{\sum_{\|x\| > t^{\sigma}} e^{c t^{\sigma (1-\epsilon)}}  f(x) Z_{x,0}^{y,t}} {\sum_{\|x\| \leq t^{\sigma}} e^{c t^{\sigma (1-\epsilon)}} f(x) Z_{x,0}^{y,t}} 
\leq \frac{\sum_{\|x\| > t^{\sigma}} e^{2c \|x\|^{1-\epsilon}} Z_{x,0}^{y,t}} {\sum_{\|x\| \leq t^{\sigma}} Z_{x,0}^{y,t}}.  
\end{equation*}
Let $\zeta, \delta > 0$. By Lemma~\ref{lm:relevant_denominator}, there are $u, T_1 > 0$ such that for every $t \geq T_1$, 
$$ 
Q \Biggl( \sum_{\|x\| \leq t^{\sigma}} Z_{x,0}^{y,t} < e^{-u} \Biggr) < \frac{\delta}{2}. 
$$ 
And by Lemma~\ref{lm:sum_of_hZ}, there is $T \geq T_1$ such that for all $t \geq T$, 
$$ 
\left \langle \Bigg \lvert \sum_{\|x\| > t^{\sigma}} e^{2c \|x\|^{1-\epsilon}} Z_{x,0}^{y,t} \Bigg \rvert \right \rangle < \frac{\delta \zeta e^{-u}}{2}. 
$$ 
For $t \geq T$, we have therefore 
\begin{align*}
&Q \left(\frac{\sum_{\|x\| > t^{\sigma}} e^{2c \|x\|^{1-\epsilon}} Z_{x,0}^{y,t}}{\sum_{\|x\| \leq t^{\sigma}} Z_{x,0}^{y,t}} > \zeta \right) \\
	\leq& Q \Biggl( \sum_{\|x\| \leq t^{\sigma}} Z_{x,0}^{y,t} < e^{-u} \Biggr) + Q \Biggl(\sum_{\|x\|>t^{\sigma}} e^{2c \|x\|^{1-\epsilon}} Z_{x,0}^{y,t} > \zeta e^{-u} \Biggr) \\
	<& \frac{\delta}{2} + \frac{e^u}{\zeta} \left \langle \Biggl \lvert \sum_{\|x\| > t^{\sigma}} e^{2c \|x\|^{1-\epsilon}} Z_{x,0}^{y,t} \Biggr\lvert \right \rangle < \delta. 
\end{align*}
This completes the proof of~\eqref{eq: convergence_C} and thus of Lemma~\ref{lm:BC_convergence}.
\epf 

\bigskip

\section{Proof of Lemma~\ref{lm:p_ratio_uniform_bound} -- Uniform bound for the quotient of transition probabilities} \label{sec:p_ratio_uniform_bound}

Lemma~\ref{lm:p_ratio_uniform_bound} is a statement about SSRW. We have no doubt that this result is known. However, since we did not find a precise reference, we decided to present a proof here. 
Fix $\sigma \in (0,1)$ and $y_1, y_2 \in \Z^d$. For an arbitrary constant $\nu \in (\tfrac{1}{2},1)$ and $t > 0$, let  
$$
J(t) := \biggl\{n \in \N: \biggl \lvert \frac{n}{t} - 1 \biggr\lvert < 1 - \nu \biggr\} = \{ n \in \N: \nu t < n < (2 - \nu)t \}. 
$$
For $t > 0$, $y \in \Z^d$, and $x \in \Z^d$ such that $\|x\| \leq t^{\sigma}$, let 
\begin{equation}\label{eq:def_Dc}
\Dc_t(y,x) := \frac{\sum_{n \notin J(t)} e^{-t} \frac{t^n}{n!} q_n^{y-x}}{\sum_{n \in J(t)} e^{-t} \frac{t^n}{n!} q_n^{y-x}}. 
\end{equation}
In~\cite[Lemma~3.5]{HKN}, we established the following result. 


\begin{lemma}     \label{lm:uniform_D_conv}
For every $y \in \Z^d$, one has  
$$ 
\lim_{t \to \infty} \sup_{\|x\| \leq t^{\sigma}} \Dc_t(y,x) = 0. 
$$ 
\end{lemma}

Recall that $q_j^z = \PP(\gamma_j = z \vert \gamma_0  = 0)$ is the transition probability for a simple symmetric random walk in discrete time on $\Z^d$. 
As shown in the proof of Lemma~3.4 in \cite{HKNN},
there are $\rho, c > 0$ such that for any $n, n' \in \N$ with $n' \leq n$ and for any $z, z' \in \Z^d$ with $\|z\| \leq \rho n$ and $\|z\|_1 \equiv n \bmod 2$, 
\begin{equation} \label{eq:Nasarov_quotient}
\frac{q_{n'}^{z'}}{q_n^z} \leq \left(1+O(n^{-\frac{2}{5}}) \right) \exp \left(c \left(\frac{\|z\|}{n} \|z-z'\| + \ln(n) \frac{n-n'}{n} \right) \right). 
\end{equation} 
Fix $\delta > 0$ and choose the parameter $\nu \in (\tfrac{1}{2},1)$ in the definition of $J(t)$ so close to $1$ that 
$$ 
\nu^{-1} < (1+\delta)^{\frac{1}{3}}. 
$$ 
With the help of Lemma~\ref{lm:uniform_D_conv}, we can choose $\tau > 0$ so large that  
\begin{align*} 
\nu^{-1} + \tau^{-1} &< (1+\delta)^{\frac{1}{3}}, \\
	1 + \sup_{\|z\| \leq t^{\sigma}} \Dc_t (y_1,z) &< (1+\delta)^{\frac{1}{3}},
 \quad \text{ for all } t \geq \tau, 
\end{align*}
and such that for all integers $n > \nu \tau$, we have $\|y_2\| + (\nu^{-1} n)^{\sigma} < \rho n$ and 
\begin{equation}    \label{eq:q_ratio_eps}  
\left(1 + O(n^{-\frac{2}{5}}) \right) \exp \left(c \left(\frac{\|y_2\| + (\nu^{-1} n)^{\sigma}}{n} \|y_1 - y_2\| + \frac{\ln(n)}{n} \right) \right) < (1+\delta)^{\frac{1}{3}}.  
\end{equation} 
For $y \in \Z^d$ and $n \in \N_0$, let 
$$ 
\iota(y,n) := \begin{cases}
                    n, & \quad \|y\|_1 \equiv n \bmod 2, \\
                    n+1, & \quad \|y \|_1 \not \equiv n \bmod 2. 
\end{cases}
$$ 
Let $t \geq \tau$ and let $x \in \Z^d$ such that $\|x\| \leq t^{\sigma}$.  For $n > \nu t$,  
$$ 
\|y_2-x\| \leq \|y_2\| + \|x\| \leq \|y_2\| + (\nu^{-1} n)^{\sigma} < \rho n \leq \rho \iota(y_2 - x,n). 
$$ 
By definition of $\iota$, we also have $\|y_2 - x\|_1 \equiv \iota(y_2 - x,n) \bmod 2$. Hence, using~\eqref{eq:Nasarov_quotient} and~\eqref{eq:q_ratio_eps}, 
\begin{align} \label{eq:quotient_q_iota}
\frac{q_n^{y_1-x}}{q_{\iota(y_2 - x,n)}^{y_2-x}} \leq& \left(1 + O(n^{-\frac{2}{5}}) \right) \exp \hspace{-1mm} \left(\hspace{-0.5mm} c \left(\frac{\|y_2-x\|}{n} \|y_1 - y_2\| + \frac{\ln(n)}{n} \right) \hspace{-0.5mm} \right) \hspace{-1mm} \\
<& \ (1+\delta)^{\frac{1}{3}}. \notag
\end{align}
Notice that for $n \in J(t)$, one has $n < (2-\nu) t \leq \nu^{-1} t$ and thus 
$$ 
\frac{t^n}{n!} 
\leq (\nu^{-1} + t^{-1}) \frac{t^{\iota(y_2 - x,n)}}{\iota(y_2 - x,n)!} 
< (1+\delta)^{\frac{1}{3}} \frac{t^{\iota(y_2 - x,n)}}{\iota(y_2 - x,n)!}.
$$ 
Then 

\begin{align}    \label{eq:p_ratio_estimate_num}  
p_t^{y_1 - x} =& (1 + \Dc_t(y_1,x)) \sum_{n \in J(t)} e^{-t} \frac{t^n}{n!} q_n^{y_1 -x}  \\
\leq& (1+\sup_{\|z\| \leq t^{\sigma}} \Dc_t(y_1,z)) e^{-t} \sum_{n \in J(t)} \frac{t^n}{n!} q_n^{y_1-x} \notag\\
\leq& (1+\delta)^{\frac{1}{3}} e^{-t} \sum_{n \in J(t): q_n^{y_1 - x} > 0} \frac{t^n}{n!} q_n^{y_1 - x} \notag \\
\leq & (1+\delta)^{\frac{2}{3}} e^{-t}\sum_{n \in J(t): q_n^{y_1 - x} > 0} \frac{t^{\iota(y_2 - x,n)}}{\iota(y_2 - x,n)!} q_{\iota(y_2 - x,n)}^{y_2-x} \frac{q_n^{y_1-x}}{q_{\iota(y_2 - x,n)}^{y_2-x}} \notag \\
\leq & (1+\delta) e^{-t}\sum_{n \in J(t): q_n^{y_1 - x} > 0} \frac{t^{\iota(y_2 - x,n)}}{\iota(y_2 - x,n)!} q_{\iota(y_2 - x,n)}^{y_2-x}, \notag
\end{align}
where in the last line we used the estimate in \eqref{eq:quotient_q_iota}. Furthermore,
\begin{equation}     \label{eq:p_ratio_estimate_denom} 
p_t^{y_2 - x}
= \sum_{n=0}^{\infty} e^{-t} \frac{t^n}{n!} q_n^{y_2-x}  
\geq e^{-t} \sum_{n \in J(t): q_n^{y_1-x} > 0} \frac{t^{\iota(y_2 - x,n)}}{\iota(y_2 - x,n)!} q_{\iota(y_2 - x,n)}^{y_2 - x}. 
\end{equation}
Combining the estimates \eqref{eq:p_ratio_estimate_num}  and \eqref{eq:p_ratio_estimate_denom} we have
$$
\frac{p_t^{y_1-x}}{p_t^{y_2-x}} < 1 + \delta.
$$
Since our choice of $\tau$ did not depend on $x$ and since $\delta > 0$ is arbitrary, this shows that 
$$ 
\limsup_{t \to \infty}  \sup_{\|x\| \leq t^{\sigma}} \frac{p_t^{y_1-x}}{p_t^{y_2-x}} \leq 1. 
$$ 
Interchanging the roles of $y_1$ and $y_2$, we obtain 
$$ 
\limsup_{t \to \infty} \sup_{\|x\| \leq t^{\sigma}} \frac{p_t^{y_2-x}}{p_t^{y_1-x}} \leq 1 
$$ 
and thus the desired result. \\

\section{Unique stationary distribution}     \label{ssec:unique_sd} 

Besides the attraction result Theorem~\ref{thm:pushforward_attraction}, a second major step towards Theorem~\ref{thm:unique_global_sol} is to establish uniqueness of the stationary distribution for a certain Markov semigroup associated with~\eqref{eq:sSHE}. 

Recall the sets $\Lc_{c, \epsilon}$, defined for any $c > 0$ and $\epsilon \in (0,1)$, consisting of functions $f$ that satisfy~\eqref{eq:sublinear_log}. Define 
\begin{equation}    \label{eq:defin_La} 
\Lc := \bigcup_{\substack{c > 0, \\ \epsilon \in (0,1)}} \Lc_{c,\epsilon},
\end{equation}
which is the set of functions $f: \Z^d \to (0,\infty)$ of subexponential growth, normalized by imposing $f(0) = 1$. Note that for any global stationary solution $Z$ to~\eqref{eq:sSHE} satisfying the subexponential growth condition of Theorem~\ref{thm:unique_global_sol}, $Q$-almost surely the quotient $Z(y,t,\omega) / Z(0,t,\omega)$ is an element of $\Lc$ for every $t \in \R$. 
We endow $\Lc$ with the subspace topology inherited from the product topology on $(0, \infty)^{\Z^d}$ and denote the corresponding Borel $\sigma$-field by $\Bc(\Lc)$. Equipped with this topology, $\Lc$ is metrizable and its topology is for instance induced by the metric 
\begin{equation}    \label{eq:def_metric_d} 
d(f,g) = \sum_{x \in \Z^d} e^{-\|x\|} \frac{\lvert f(x) - g(x) \rvert}{1 + \lvert f(x) - g(x) \rvert}. 
\end{equation} 
The metric space $(0, \infty)^{\Z^d}$ with the metric $d$ is separable, and so is $\Lc$ as a topological subspace of $(0,\infty)^{\Z^d}$. For $\omega \in \Omega$, $s, t \in \R$ such that $s \leq t$, $f \in \Lc$, and $y \in \Z^d$, define 
$$
L^{s,t}_{\omega} f(y) := \frac{\sum_{x \in \Z^d} f(x) Z_{x,s}^{y,t}(\omega)}{\sum_{x \in \Z^d} f(x) Z_{x,s}^{0,t}(\omega)}   
$$
whenever the expression on the right-hand side is well-defined -- this is in particular the case if the term in the denominator is finite. 

The following three lemmas are proved in Section~\ref{sec:skew_product_proofs}.

\begin{lemma}   \label{lm:L_invariance} 
The set $\Lc$ is $Q$-almost surely invariant under the dynamics induced by $L$, i.e., $Q$-almost surely for every $f \in \Lc$ and for every $s, t \in \R$ such that $s \leq t$, we have $L^{s,t}_{\omega} f \in \Lc$.  
\end{lemma}

Recall the set $\Omega^{\textup{sol}}$, whose existence was postulated in Proposition~\ref{prop:global_stationary_sol}.

\begin{lemma}    \label{lm:Omega_subset} 
If $\beta$ is sufficiently small, there is an $\Fc$-measurable set $\widetilde{\Omega} \subset \Omega^{\textup{sol}}$ with $Q(\widetilde{\Omega}) = 1$ that satisfies the following conditions:  
\begin{enumerate}
\item $\widetilde{\Omega}$ is invariant under $\theta_s$ for every $s \in \R$, i.e., $\theta_s(\widetilde{\Omega}) = \widetilde{\Omega}$ for every $s \in \R$; 
\item For every $\omega \in \widetilde{\Omega}$, $f \in \Lc$, and $s, t \in \R$ such that $s \leq t$, one has $L^{s,t}_{\omega} f \in \Lc$;
\item For every $\omega \in \widetilde{\Omega}$, 
$$
\lim_{s \to -\infty} Z_s^{y,0}(\omega) = Z_{-\infty}^{y,0}(\omega) > 0, \qquad \forall y \in \Z^d;
$$
\item For every $\omega \in \widetilde{\Omega}$, the function 
$$
\mathcal{Y}(\omega): \Z^d \to (0, \infty), \ y \mapsto Z_{-\infty}^{y,0}(\omega)/Z_{-\infty}^{0,0}(\omega)
$$
is an element of $\Lc$. 
\end{enumerate}
\end{lemma}  

\begin{lemma} \label{lm:cocycle} 
The map $\varphi: [0, \infty) \times \widetilde{\Omega} \times \Lc \to \Lc$, given by 
$$
(t,\omega,f) \mapsto \varphi^t_{\omega} f := L^{0,t}_{\omega} f, 
$$
defines a cocycle; i.e., 
\begin{enumerate}
    \item $\varphi$ is $(\Bc([0,\infty)) \otimes \widetilde{\Fc}_{0,\infty} \otimes \Bc(\Lc), \Bc(\Lc))$-measurable, where $\Bc([0, \infty))$ is the Borel $\sigma$-field on $[0,\infty)$ and $\widetilde{\Fc}_{0,\infty}$ is the restriction of $\sigma(W^x_u: u \geq 0, x \in \Z^d)$ to $\widetilde{\Omega}$; 
    \item For all $\omega \in \widetilde{\Omega}$, $\varphi_{\omega}^0$ is the identity on $\Lc$ and 
    $$
    \varphi_{\omega}^{s+t} = \varphi^t_{\theta_s \omega} \circ \varphi^s_{\omega}, \quad \forall s, t \geq 0. 
    $$
\end{enumerate}
\end{lemma} 

On $\Lc$, consider the Markov semigroup $(\Pp^t)_{t \geq 0}$ defined by 
\begin{equation}    \label{eq:Markov_semi} 
\Pp^t(f,B) := Q \left(\left\{\omega \in \widetilde{\Omega}: \ \varphi^t_{\omega} f \in B \right\} \right), \qquad f \in \Lc, \ B \in \Bc(\Lc). 
\end{equation}

\begin{definition}   \rm 
    A stationary distribution for $(\Pp^t)_{t \geq 0}$ is a probability measure $\nu$ on $(\Lc, \Bc(\Lc))$ which satisfies 
    $$
\nu \Pp^t(\cdot) := \int_{\Lc} \Pp^t(f, \cdot) \ \nu(df) = \nu(\cdot), \qquad \forall t \geq 0. 
$$
\end{definition}

The following theorem asserts existence and uniqueness for the stationary distribution associated with $(\Pp^t)_{t \geq 0}$. We prove it in Section~\ref{sec:proof_unique_sd}.

\begin{theorem}    \label{thm:unique_stationary_dis} 
If $\beta$ is sufficiently small, the Markov semigroup $(\Pp^t)_{t \geq 0}$ admits a unique stationary distribution $\nu^{\star}$ given by 
$$
\nu^{\star}(\cdot) = \int_{\widetilde{\Omega}} \delta_{\mathcal{Y}(\omega)}(\cdot) \ Q(d \omega),  
$$
where $\mathcal{Y}(\omega)$ was introduced in Lemma~\ref{lm:Omega_subset} and where $\delta_{\mathcal{Y}(\omega)}$ is the Dirac measure concentrated on the function $\mathcal{Y}(\omega)$. 
\end{theorem}





\section{Proof of Theorem~\ref{thm:unique_global_sol}}  \label{sec:proof_unique_global_sol}

Let $Z: \Z^d \times \R \times \Omega' \to (0,\infty)$ be a global stationary solution to~\eqref{eq:sSHE} such that for every $\omega \in \Omega'$ there is $\epsilon \in (0,1)$ with 
$$
\lim_{\|y\| \to \infty} \frac{\lvert \ln(Z(y,0,\omega)) \rvert}{\|y\|^{1-\epsilon}} = 0. 
$$
For $\omega \in \Omega'$, let $\mathcal{Z}(\omega)$ denote the function 
$$
\mathcal{Z}(\omega): \Z^d \to (0,\infty), \ y \mapsto Z(y,0,\omega)/Z(0,0,\omega). 
$$
Then $\mathcal{Z}(\omega) \in \Lc$, where $\Lc$ was defined in~\eqref{eq:defin_La}, Section~\ref{ssec:unique_sd}. Recall that in Lemma~\ref{lm:Omega_subset} we defined the function 
$$
\mathcal{Y}(\omega): \Z^d \to (0, \infty), \ y \mapsto Z_{-\infty}^{y,0}(\omega)/Z_{-\infty}^{0,0}(\omega)
$$
for $\omega \in \widetilde{\Omega}$. According to Theorem~\ref{thm:unique_stationary_dis}, for $\beta$ sufficiently small, 
$$
\nu^{\star}(\cdot) = \int_{\widetilde{\Omega}} \delta_{\mathcal{Y}(\omega)}(\cdot) \ Q(d \omega) 
$$
is the unique stationary distribution for the Markov semigroup $(\Pp^t)_{t \geq 0}$ on $(\Lc, \Bc(\Lc))$ defined in~\eqref{eq:Markov_semi}, Section~\ref{ssec:unique_sd}. 

As in the proof of Theorem~\ref{thm:unique_stationary_dis}, one shows that 
$$
\nu(\cdot) := \int_{\Omega'} \delta_{\mathcal{Z}(\omega)}(\cdot) \ Q(d \omega) 
$$
is a stationary distribution for $(\Pp^t)_{t \geq 0}$, so $\nu = \nu^{\star}$ by the uniqueness part of Theorem~\ref{thm:unique_stationary_dis}. As in the proof of the Ledrappier--Le Jan--Crauel correspondence theorem (see Theorem~4.2.9 in~\cite{Kuksin_Shirikyan}, and also~\cite{ledrappier2006positivity} for the original reference), one proceeds to show that $\delta_{\mathcal{Z}(\omega)} = \delta_{\mathcal{Y}(\omega)}$ and hence $\mathcal{Z}(\omega) = \mathcal{Y}(\omega)$ for $Q$-almost every $\omega$. 
We are not quite done yet because $\mathcal{Z}$ and $\mathcal{Y}$ do not depend on $t$. To complete the proof, let $\Omega_1 \subset \widetilde{\Omega} \cap \Omega'$ such that $Q(\Omega_1) = 1$ and $\mathcal{Z}(\omega) = \mathcal{Y}(\omega)$ for every $\omega \in \Omega_1$. To deal with negative times $t$, set 
$$
\Omega_2 := \bigcap_{r \in \N_0} \theta_r(\Omega_1). 
$$
Since $Q$ is invariant under each shift $\theta_r$, one has $Q(\Omega_2) = 1$. Moreover, for every $u \in \N_0$, 
$$
\theta_{-u}(\Omega_2) = \bigcap_{r \in \N_0} \theta_{-u+r}(\Omega_1) \subset \Omega_2. 
$$
Finally, as $\Omega_2 \subset \Omega_1$, one has $\mathcal{Z}(\omega) = \mathcal{Y}(\omega)$ for every $\omega \in \Omega_2$. 

We claim that for every $\omega \in \Omega_2$ and for all $(y,t) \in \Z^d \times \R$, 
$$
Z(y,t,\omega) = Z_{-\infty}^{y,t}(\omega) \cdot \frac{Z(0,0,\omega)}{Z_{-\infty}^{0,0}(\omega)}. 
$$
We first show that 
$$
Z(y,t,\omega) = Z_{-\infty}^{y,t}(\omega) \cdot \frac{Z(0,t,\omega)}{Z_{-\infty}^{0,t}(\omega)}, \quad \forall \omega \in \Omega_2, \ \forall (y,t) \in \Z^d \times \R, 
$$
and then prove that the factor $Z(0,t,\omega)/Z_{-\infty}^{0,t}(\omega)$ is actually constant in $t$. Since both $Z$ and $Z_{-\infty}^{\cdot, \cdot}(\cdot)$ are global stationary solutions to~\eqref{eq:sSHE}, we have for every $\omega \in \Omega_2$, $y \in \Z^d$, and $t > 0$ 
\begin{align*}
\frac{Z(y,t,\omega)}{Z(0,t,\omega)} =& \frac{\sum_{x \in \Z^d} Z(x,0,\omega) Z_{x,0}^{y,t}(\omega)}{\sum_{x \in \Z^d} Z(x,0,\omega) Z_{x,0}^{0,t}(\omega)} \\
=& \frac{\sum_{x \in \Z^d} \mathcal{Z}(\omega)[x] Z_{x,0}^{y,t}(\omega)}{\sum_{x \in \Z^d} \mathcal{Z}(\omega)[x] Z_{x,0}^{0,t}(\omega)} = \frac{\sum_{x \in \Z^d} \mathcal{Y}(\omega)[x] Z_{x,0}^{y,t}(\omega)}{\sum_{x \in \Z^d} \mathcal{Y}(\omega)[x] Z_{x,0}^{0,t}(\omega)} = \frac{Z^{y,t}_{-\infty}(\omega)}{Z^{0,t}_{-\infty}(\omega)}. 
\end{align*} 
For any $t \leq 0$ there is $u \in \N_0$ such that $t + u > 0$, so we obtain for every $\omega \in \Omega_2$ and $y \in \Z^d$ 
$$
\frac{Z(y,t,\omega)}{Z(0,t,\omega)} = \frac{Z(y,t+u,\theta_{-u} \omega)}{Z(0,t+u, \theta_{-u} \omega)} = \frac{Z^{y,t+u}_{-\infty}(\theta_{-u} \omega)}{Z^{0,t+u}_{-\infty}(\theta_{-u} \omega)} = \frac{Z^{y,t}_{-\infty}(\omega)}{Z^{0,t}_{-\infty}(\omega)}. 
$$
It remains to prove that $Z(0,t,\omega)/Z_{-\infty}^{0,t}(\omega)$ is constant in $t$. Let $\omega \in \Omega_2$ and let $s, t \in \R$ such that $s < t$. Choose $u \in \N_0$ so large that $-u < s$. Then  
\begin{align*}
\frac{Z(0,s,\omega)}{Z(0,-u,\omega)} =& \sum_{x \in \Z^d} \frac{Z(x,-u,\omega)}{Z(0,-u,\omega)} Z_{x,-u}^{0,s}(\omega) = \sum_{x \in \Z^d} \frac{Z(x,0,\theta_{-u} \omega)}{Z(0,0, \theta_{-u} \omega)} Z_{x,-u}^{0,s}(\omega) \\
=& \sum_{x \in \Z^d} \mathcal{Z}(\theta_{-u} \omega)[x] Z_{x,-u}^{0,s}(\omega) = \sum_{x \in \Z^d} \mathcal{Y}(\theta_{-u} \omega)[x] Z_{x,-u}^{0,s}(\omega) = \frac{Z_{-\infty}^{0,s}(\omega)}{Z_{-\infty}^{0,-u}(\omega)}. 
\end{align*} 
By the same argument, we also have 
$$
\frac{Z(0,t,\omega)}{Z(0,-u,\omega)} = \frac{Z_{-\infty}^{0,t}(\omega)}{Z_{-\infty}^{0,-u}(\omega)}, 
$$
so $Z(0,s,\omega)/Z_{-\infty}^{0,s}(\omega) = Z(0,t,\omega)/Z_{-\infty}^{0,t}(\omega)$. This shows that $Z(0,t,\omega)/Z_{-\infty}^{0,t}(\omega)$ is constant in $t$. 

\bigskip

\section{Existence and uniqueness of the stationary distribution -- Proof of Theorem~\ref{thm:unique_stationary_dis}}    \label{sec:proof_unique_sd}

\subsection{Existence} 

Let us show that 
$$
\nu^{\star}(\cdot) = \int_{\widetilde{\Omega}} \delta_{\mathcal{Y}(\omega)}(\cdot) \ Q(d \omega) 
$$
is indeed a stationary distribution for $(\Pp^t)_{t \geq 0}$. 
Let $t \geq 0$ and $B \in \Bc(\Lc)$. One has 
\begin{align*} 
\nu^{\star} \Pp^t(B) =& \int_{\Lc} \Pp^t(f,B) \ \nu^{\star}(df) \\
=& \int_{\widetilde{\Omega}} \Pp^t(\mathcal{Y}(\omega), B) \ Q(d \omega) \\
=& \int_{\widetilde{\Omega}} Q(\{\alpha \in \widetilde{\Omega}: \ \varphi^t_{\alpha} \mathcal{Y}(\omega) \in B\}) \ Q(d \omega) \\
=& \int_{\widetilde{\Omega}} \int_{\widetilde{\Omega}} \id_B(\varphi^t_{\alpha} \mathcal{Y}(\omega)) \ Q(d \alpha) \ Q(d \omega) \\
=& \int_{\widetilde{\Omega}} \id_B(\varphi^t_{\omega} \mathcal{Y}(\omega)) \ Q(d \omega) \\
=& \int_{\widetilde{\Omega}} \id_B(\mathcal{Y}(\theta_t \omega)) \ Q(d \omega) = \nu^{\star}(B). 
\end{align*} 

\medskip 

\subsection{Uniqueness}

Let $\nu$ be any stationary distribution for $(\Pp^t)_{t \geq 0}$. We will show that, as $t \to \infty$, $\nu \Pp^t$ converges weakly to $\nu^{\star}$. This implies uniqueness because $\nu \Pp^t = \nu$ and because any weakly convergent sequence of Borel probability measures on a metric space has a unique limit. 

Let $F: \Lc \to \R$ be bounded and Lipschitz continuous, i.e., there exists $L > 0$ such that 
$$
\lvert F(f) - F(g) \rvert \leq L d(f,g), \quad \forall f, g \in \Lc, 
$$
where $d$ is the metric on $\Lc$ defined in~\eqref{eq:def_metric_d}. One has 
\begin{align*}
    \int_{\Lc} F(f) \ \nu \Pp^t(df) =& \int_{\Lc} \Pp^t F(f) \ \nu(df) \\
    =& \int_{\Lc} \int_{\Lc} F(g) \ \Pp^t(f,dg) \ \nu(df) \\
    =& \int_{\Lc} \int_{\Omega} F(\varphi^t_{\omega} f) \ Q(d \omega) \ \nu(df), 
\end{align*}
where, in the last step, we used the change-of-variables formula. In addition, 
\begin{align*}
    \int_{\Lc} F(f) \ \nu^{\star}(df) =& \int_{\Omega} F(\mathcal{Y}(\omega)) \ Q(d \omega) \\
    =& \int_{\Omega} F(\mathcal{Y}(\theta_t \omega)) \ Q(d \omega) \\
    =& \int_{\Lc} \int_{\Omega} F(Z_{-\infty}^{\cdot, t}(\omega)/Z_{-\infty}^{0,t}(\omega)) \ Q(d \omega) \ \nu(df). 
\end{align*}
As a result, 
\begin{align*}
    &\biggl\lvert \int_{\Lc} F(f) \ \nu \Pp^t(df) - \int_{\Lc} F(f) \ \nu^{\star}(df) \biggr\rvert \\
    \leq& \int_{\Lc} \int_{\Omega} \left\lvert F(\varphi^t_{\omega} f) - F(Z_{-\infty}^{\cdot, t}(\omega)/Z_{-\infty}^{0,t}(\omega)) \right\rvert \ Q(d \omega) \ \nu(df).  
    \end{align*}
Let us show that for every $f \in \Lc$
\begin{equation}    \label{eq:convergence_fixed_f} 
\lim_{t \to \infty} \int_{\Omega} \left\lvert F(\varphi^t_{\omega} f) - F(Z_{-\infty}^{\cdot,t}(\omega)/Z_{-\infty}^{0,t}(\omega)) \right\rvert \ Q(d \omega) = 0. 
\end{equation} 
Fix $f \in \Lc$ and let $\varepsilon > 0$. Let $Z \subset \Z^d$ be a finite set such that 
$$
\sum_{y \in \Z^d \setminus Z} e^{-\|y\|} < \frac{\varepsilon}{3L}. 
$$
By Theorem~\ref{thm:pushforward_attraction}, 
$$
\left\lvert \frac{u_f(y,t)}{u_f(0,t)} - \frac{Z_{-\infty}^{y,t}}{Z_{-\infty}^{0,t}} \right\rvert \xrightarrow[t \to \infty] {} 0  \text{ in probability,} \quad \forall y \in Z.  
$$
Let $\lvert Z \rvert$ denote the cardinality of $Z$. Since $Z$ is finite, there exists $T > 0$ such that 
$$
Q \biggl( \biggl\lvert \frac{u_f(y,t)}{u_f(0,t)} - \frac{Z_{-\infty}^{y,t}}{Z_{-\infty}^{0,t}} \biggr\rvert > \frac{\varepsilon}{3 \lvert Z \rvert L} \biggr) < \frac{\varepsilon}{3 \lvert Z \rvert L}
$$
for every $t \geq T$ and for every $y \in Z$. For $y \in Z$ and $t \geq T$, set 
$$
A(y,t) = \biggl\{\biggl \lvert \frac{u_f(y,t)}{u_f(0,t)} - \frac{Z_{-\infty}^{y,t}}{Z_{-\infty}^{0,t}} \biggr \rvert > \frac{\varepsilon}{3 \lvert Z \rvert L} \biggr\}. 
$$
As $\varphi^t_{\omega} f(y) = u_f(y,t)/u_f(0,t)$, one has for $t \geq T$
\begin{align*}
& \int_{\Omega} \left \lvert F(\varphi^t_{\omega} f) - F(Z_{-\infty}^{\cdot,t}(\omega)/Z_{-\infty}^{0,t}(\omega)) \right\rvert \ Q(d \omega) \\
    \leq& \ L \int_{\Omega} d \biggl(\varphi^t_{\omega} f, \frac{Z_{-\infty}^{\cdot,t}(\omega)}{Z_{-\infty}^{0,t}(\omega)} \biggr) \ Q(d \omega) \\ 
    \leq& \ L \int_{\Omega} \sum_{y \in \Z^d \setminus Z} e^{-\|y\|} \ Q(d \omega) + L \sum_{y \in Z} e^{-\|y\|} Q(A(y,t)) \\
    &+ L \sum_{y \in Z} e^{-\|y\|} \int_{\Omega \setminus A(y,t)} \biggl\lvert \frac{u_f(y,t)}{u_f(0,t)} - \frac{Z_{-\infty}^{y,t}(\omega)}{Z_{-\infty}^{0,t}(\omega)} \biggr \rvert \ Q(d \omega) \\
    <& \frac{\varepsilon}{3} + \frac{\varepsilon}{3} + \frac{\varepsilon}{3} = \varepsilon. 
\end{align*} 
We have thus proved~\eqref{eq:convergence_fixed_f}. Since $F$ is bounded, dominated convergence yields 
$$
\lim_{t \to \infty} \int_{\Lc} \int_{\Omega} \left \lvert F(\varphi^t_{\omega} f) - F(Z_{-\infty}^{\cdot,t}(\omega)/Z_{-\infty}^{0,t}(\omega)) \right \rvert \ Q(d \omega) \ \nu(df) = 0. 
$$
Therefore, 
$$
\lim_{t \to \infty} \int_{\Lc} F(f) \ \nu \Pp^t(df) = \int_{\Lc} F(f) \ \nu^{\star}(df). 
$$
As this convergence holds for every bounded Lipschitz continuous function on $\Lc$, the Portmanteau theorem (see, e.g., Theorem~4.1 in~\cite{benaim2022markov}) implies that $\nu \Pp^t$ converges weakly to $\nu^{\star}$. 

\bigskip

\section{Proofs of Lemmas~\ref{lm:L_invariance},~\ref{lm:Omega_subset}, and~\ref{lm:cocycle}}  \label{sec:skew_product_proofs} 

\subsection{Proof of Lemma~\ref{lm:L_invariance}}

We begin by defining some auxiliary sets. Consider the function 
$$
g: \Z^d \times \Z^d \to (0, \infty), \ (x,y) \mapsto (1+\|x\|)^{d+1} (1+\|y\|)^{d+1} 
$$
and notice that 
$$
\sum_{x, y \in \Z^d} \frac{1}{g(x,y)} < \infty. 
$$
For $M \in \N$, let 
\begin{align*}
\Omega^1_M :=& \biggl\{\omega \in \Omega: \exists R_M(\omega) > 0 \ \text{s.t.} \ \forall x, y \in \Z^d \ \text{with} \ \max\{\|x\|, \|y\|\} > R_M(\omega): \\
& \sup_{\substack{s, t \in (-M,M), \\ s < t}} Z_{x,s}^{y,t}(\omega) \leq (p_{2M}^{y-x})^{\frac{1}{2}} g(x,y), \ \inf_{\substack{s,t \in (-M,M), \\ s < t}} Z_{x,s}^{x,t}(\omega) \geq g(x,x)^{-1} \biggr\}. 
\end{align*}
Finally, set $\Omega^1 := \bigcap_{M \in \N} \Omega^1_M$. We first show that $Q(\Omega^1_M) = 1$ for every $M \in \N$, and hence $Q(\Omega^1) = 1$. Fix $M \in \N$. For $x, y \in \Z^d$, consider the events 
\begin{align*}
E_1(x,y) :=& \left\{\omega \in \Omega: \sup_{s, t \in (-M,M), s < t} Z_{x,s}^{y,t}(\omega) > (p_{2M}^{y-x})^{\frac{1}{2}} g(x,y) \right\}, \\
E_2(x) :=& \left\{\omega \in \Omega: \inf_{s, t \in (-M,M), s < t} Z_{x,s}^{x,t}(\omega) < g(x,x)^{-1} \right\}. 
\end{align*}
By the first Borel--Cantelli lemma, $Q(\Omega^1_M) = 1$ will follow once we show that 
$$
\sum_{x, y \in \Z^d} Q(E_1(x,y)) + \sum_{x \in \Z^d} Q(E_2(x)) < \infty. 
$$
By Markov's inequality, 
$$
Q(E_1(x,y)) \leq \frac{(p_{2M}^{y-x})^{-\frac{1}{2}} \langle \sup_{s, t \in (-M,M), s < t} Z_{x,s}^{y,t} \rangle}{g(x,y)}. 
$$
We will now show that the expression in the numerator is bounded by a constant that only depends on $M$ and $\beta$. For fixed $\omega \in \Omega$ and $s, t \in (-M,M)$ such that $s < t$, we can estimate 
\begin{align}   \label{eq:f_E_estim}
Z_{x,s}^{y,t}(\omega) <& p_{t-s}^{y-x} \E_{x,s}^{y,t} e^{\beta \Ac_s^t(\cdot, \omega)} \\
=&\frac{1}{p_{s+M}^0} \int \id_{\eta_s = x, \eta_t = y} e^{\beta \Ac_s^t(\eta, \omega)} \ \Pp_{x,-M}(d \eta) \notag \\
\leq& \left(\frac{p_{t-s}^{y-x}}{p_{s+M}^0} \right)^{\frac{1}{2}} \left(\int e^{2 \beta \lvert \Ac_s^t(\eta, \omega) \rvert} \ \Pp_{x,-M}(d \eta) \right)^{\frac{1}{2}}, \notag
\end{align}
where the integral is taken over possible realizations of $\eta$ and where we used the trivial estimate $e^{-\frac{\beta^2}{2} (t-s)} < 1$ to obtain the first inequality. For every $s \in (-M,M)$, we have 
\begin{equation}  \label{eq:p_s_M_lower_bound} 
p_{s+M}^0 \geq e^{-(s+M)} \geq e^{-2M}. 
\end{equation} 
Moreover, for every $z \in \Z^d \setminus \{0\}$, 
$$
\frac{d}{dr} p_r^z \geq 0, \quad \forall r \in (0, \|z\|_1). 
$$
To see this, recall that for every $z \in \Z^d$ and $r \geq 0$, 
$$
p_r^z = e^{-r} \sum_{n=0}^{\infty} \frac{r^n}{n!} q_n^z. 
$$
Differentiating the expression on the right with respect to $r$ yields 
$$
e^{-r} \left(-q_0^z + \sum_{n=1}^{\infty} q_n^z \left(\frac{r^{n-1}}{(n-1)!} - \frac{r^n}{n!} \right) \right). 
$$
Since $z \neq 0$, we have $q^z_0 = 0$. For every $n \geq r$, 
$$
\frac{r^n}{n!} = \frac{r}{n} \cdot \frac{r^{n-1}}{(n-1)!} \leq \frac{r^{n-1}}{(n-1)!}. 
$$
And for $n < r < \|z\|_1$, we have $q_n^z = 0$. This proves the claim. For every $x, y \in \Z^d$ such that $\|y-x\|_1 > 2M$ and $s, t \in (-M,M)$ with $s < t$, we have thus 
$$
p^{y-x}_{t-s} \leq p^{y-x}_{2M}. 
$$
Hence, for every $x, y \in \Z^d$ and for every $s, t \in (-M,M)$ with $s < t$, 
$$
p^{y-x}_{t-s} \leq \id_{\|y-x\|_1 > 2M} p^{y-x}_{2M} + \id_{\|y-x\|_1 \leq 2M}.  
$$
Together with~\eqref{eq:f_E_estim} and~\eqref{eq:p_s_M_lower_bound}, this yields 
\begin{align} \label{eq:p_Z_estim} 
&(p_{2M}^{y-x})^{-\frac{1}{2}} Z_{x,s}^{y,t}(\omega) \\
\leq& e^M  (1 + \max_{z \in \Z^d: \|z\|_1 \leq 2M} (p_{2M}^z)^{-\frac{1}{2}}) \left(\int e^{2 \beta \lvert \Ac_s^t(\eta, \omega) \rvert} \ \Pp_{x,-M}(d \eta) \right)^{\frac{1}{2}}.  \notag
\end{align} 
To make the rest of the exposition easier to read, we write $n(M)$ instead of $n_{-M,M}$. Let $\eta$ be a continuous-time random-walk path starting from $x$ at time $-M$, and observed on the time interval $[-M,M]$, with jump sites $\gamma'_1, \ldots, \gamma'_{n(M)}$, and jump times $s'_1 < \ldots < s'_{n(M)}$. Assume in addition that $\eta_s = x$ and $\eta_t = y$. Set $\gamma'_0 := x$, $s'_0 := -M$, and $s'_{n(M) +1} := M$. We use the following relabelling of the jump times and jump sites: Let $j^* := \min\{j: s'_j > s\}$ and 
$$
s_k := s'_{j^* + k-1}, \quad \gamma_k := \gamma'_{j^* + k-1}, \quad \forall k \in \{1, \ldots, n_{s,t}\}. 
$$
Also set $\gamma_0 := x$, $s_0 := s$, and $s_{n_{s,t} +1} := t$. Then 
\begin{align*} 
\lvert \Ac_s^t(\eta, \omega) \rvert \leq& \sum_{j=0}^{n_{s,t}} \left \lvert \omega(\gamma_j, s_{j+1}) - \omega(\gamma_j, s_j) \right \rvert \\ 
\leq& 4 \max_{z \in \{x,y\}} \max_{r \in [-M,M]} \left \lvert \omega(z, r) \right \rvert + \sum_{j=0}^{n(M)} \left \lvert \omega(\gamma'_j, s'_{j+1}) - \omega(\gamma'_j, s'_j) \right \rvert.  
\end{align*}
Since the expression on the right does not depend on $s$ or $t$, we obtain 
\begin{align*}
& \left \langle \sup_{s, t \in (-M,M), s < t} \left( \int e^{2 \beta \lvert \Ac_s^t(\eta, \cdot) \rvert} \Pp_{x,-M}(d \eta) \right)^{\frac{1}{2}} \right \rangle \\
\leq& \biggl \langle \biggl( \int \exp \biggl(2 \beta \biggl(4 \max_{z \in \{x,y\}} \max_{r \in [-M,M]} \left \lvert \omega(z, r) \right \rvert \\
&+ \sum_{j=0}^{n(M)} \left \lvert \omega(\gamma'_j, s'_{j+1}) - \omega(\gamma'_j, s'_j) \right \rvert \biggr) \biggr) \Pp_{x,-M}(d \eta) \biggr)^{\frac{1}{2}} \biggr \rangle. 
\end{align*}
By Jensen's inequality, Fubini, and Cauchy--Schwarz, the expression on the right is bounded from above by
\begin{align}   \label{eq:Jen_Fub_CS} 
\biggl( \int & \biggl \langle \exp \biggl(16 \beta \max_{z \in \{x,y\}}  \max_{r \in [-M,M]} \left \lvert \omega(z,r) \right \rvert \biggr) \biggr\rangle^{\frac{1}{2}} \\
& \biggl \langle \exp \biggl(4 \beta \sum_{j=0}^{n(M)} \left \lvert \omega(\gamma'_j, s'_{j+1}) - \omega(\gamma'_j, s'_j) \right \rvert \biggr) \biggr \rangle^{\frac{1}{2}} \ \Pp_{x,-M}(d \eta) \biggr)^{\frac{1}{2}}.  \notag 
\end{align}
For a fixed realization $\eta$ of the random walk and for $0 \leq j \leq n(M)$, we have 
$$
\left \langle \exp \left(4 \beta \left \lvert \omega(\gamma'_j, s'_{j+1}) - \omega(\gamma'_j, s'_j) \right \rvert \right) \right \rangle < 2 e^{8 (s'_{j+1}-s'_j) \beta^2}. 
$$
Thus, 
$$
\biggl \langle \exp \biggl(4 \beta \sum_{j=0}^{n(M)} \left \lvert \omega(\gamma'_j, s'_{j+1}) - \omega(\gamma'_j, s'_j) \right \rvert \biggr) \biggr \rangle \leq 2^{n(M)+1} e^{16 M \beta^2}. 
$$
For $z \in \{x, y\}$, let $\xi_z := \sup_{r \in [-M,M]} \lvert \omega(z,r) \rvert$, and 
\begin{align*}
\xi_z^+(1) :=& \max_{r \in [0,M]} \omega(z,r), & \xi_z^+(2) :=& \max_{r \in [0,M]} (- \omega(z,r)), \\
\xi_z^-(1) :=& \max_{r \in [-M,0]} \omega(z,r), & \xi_z^-(2) :=& \max_{r \in [-M,0]} (-\omega(z,r)). 
\end{align*}
Then
$$
\max_{z \in \{x,y\}} \max_{r \in [-M,M]} \left \lvert \omega(z,r) \right \rvert \leq \xi_x + \xi_y \leq \sum_{z \in \{x,y\}} \sum_{i \in \{1,2\}} (\xi_z^+(i) + \xi_z^-(i)).  
$$
It follows that 
\begin{align*}
&\biggl \langle \exp \biggl(16 \beta \max_{z \in \{x,y\}} \max_{r \in [-M,M]} \left \lvert \omega(z,r) \right \rvert \biggr) \biggr \rangle \\
\leq& \biggl \langle \exp \biggl(16 \beta \sum_{z \in \{x,y\}} \sum_{i \in \{1,2\}} (\xi_z^+(i) + \xi_z^-(i)) \biggr) \biggr \rangle \\
=& \prod_{z \in \{x,y\}} \left \langle e^{16 \beta \xi_z^+(1)} e^{16 \beta \xi_z^+(2)} \right \rangle \left \langle e^{16 \beta \xi_z^-(1)} e^{16 \beta \xi_z^-(2)} \right \rangle \\
\leq& \prod_{z \in \{x,y\}} \prod_{i \in \{1,2\}} \left \langle e^{32 \beta \xi_z^+(i)} \right \rangle^{\frac{1}{2}} \left \langle e^{32 \beta \xi_z^-(i)} \right \rangle^{\frac{1}{2}} = \left \langle e^{32 \beta \xi_0^+(1)} \right \rangle^4.  
\end{align*}
The random variable $\xi_0^+(1)$ has the same distribution as absolute value of a Gaussian random variable with mean $0$ and variance $M$. Therefore, 
$$
\left \langle e^{32 \beta \xi_0^+(1)} \right \rangle \leq 2 e^{512 M \beta^2}.  
$$
We have thus shown the estimate 
\begin{align*}
\biggl \langle \exp \biggl(16 \beta & \max_{z \in \{x,y\}}  \max_{r \in [-M,M]} \left \lvert \omega(z,r) \right \rvert \biggr) \biggr \rangle^{\frac{1}{2}} \\
&\biggl \langle \exp \biggl(4 \beta \sum_{j=0}^{n(M)} \left \lvert \omega(\gamma'_j, s'_{j+1})  - \omega(\gamma'_j, s'_j) \right \rvert \biggr) \biggr \rangle^{\frac{1}{2}} 
\leq 4 e^{1032 M \beta^2} 2^{\frac{n(M)+1}{2}}.  
\end{align*}
As a result, the expression in~\eqref{eq:Jen_Fub_CS} is bounded from above by 
\begin{align*}
\biggl(\int 4 e^{1032 M \beta^2} 2^{\frac{n(M)+1}{2}} \ \Pp_{x,-M}(d \eta) \biggr)^{\frac{1}{2}} 
=& 2 e^{516 M \beta^2} \biggl(\sum_{n=0}^{\infty} e^{-2M} \frac{(2M)^n}{n!} 2^{\frac{n+1}{2}} \biggr)^{\frac{1}{2}} \\
=& 2^{\frac{5}{4}}  \cdot e^{516 M \beta^2} \exp(M (\sqrt{2}-1)) =: h(M). 
\end{align*}
Hence, together with~\eqref{eq:p_Z_estim}, 
$$
Q(E_1(x,y)) \leq \frac{H(M)}{g(x,y)}, 
$$
where 
$$
H(M) := e^M \left(1 + \max_{z \in \Z^d: \|z\|_1 \leq 2M} (p_{2M}^z)^{-\frac{1}{2}} \right) h(M). 
$$
Furthermore, by Markov's inequality and Jensen's inequality, 
\begin{align}   \label{eq:Markov_E_2}
Q(E_2(x)) \leq& \frac{\langle \sup_{s, t \in (-M,M), s < t} (1/ Z_{x,s}^{x,t}) \rangle}{g(x,x)} \\
\leq& e^{\beta^2 M} \frac{\langle \sup_{s,t \in (-M,M), s < t} (p_{t-s}^0)^{-1} \E_{x,s}^{x,t} e^{-\beta \Ac_s^t} \rangle}{g(x,x)}. \notag
\end{align}
For fixed $\omega \in \Omega$ and fixed $s,t \in (-M,M)$ such that $s < t$, we have  
\begin{align*}
\frac{1}{p_{t-s}^0} \E_{x,s}^{x,t} e^{-\beta \Ac_s^t(\cdot, \omega)} \leq& \left(p^0_{s+M} \right)^{-\frac{1}{2}} \left(p^0_{t-s} \right)^{-\frac{3}{2}} \biggl(\int e^{2 \beta \lvert \Ac_s^t(\eta, \omega) \rvert} \ \Pp_{x,-M}(d \eta) \biggr)^{\frac{1}{2}} \\
\leq& e^{4M} \biggl( \int e^{2 \beta \lvert \Ac_s^t(\eta, \omega) \rvert} \ \Pp_{x,-M}(d \eta) \biggr)^{\frac{1}{2}}, 
\end{align*}
where we used the estimate in~\eqref{eq:p_s_M_lower_bound}. We have already established that 
$$
\biggl \langle \sup_{s, t \in (-M,M), s < t} \left(\int e^{2 \beta \lvert \Ac_s^t(\eta, \cdot) \rvert} \ \Pp_{x,-M}(d \eta) \right)^{\frac{1}{2}} \biggr \rangle \leq h(M). 
$$
Hence, the expression on the right-hand side of~\eqref{eq:Markov_E_2} is bounded from above by 
$$
\frac{e^{\beta^2 M + 4M} h(M)}{g(x,x)}. 
$$
Then 
$$
\sum_{x, y \in \Z^d} Q(E_1(x,y)) \leq H(M) \sum_{x, y \in \Z^d} \frac{1}{g(x,y)} < \infty 
$$
and 
$$
\sum_{x \in \Z^d} Q(E_2(x)) \leq e^{\beta^2 M + 4M} h(M) \sum_{x \in \Z^d} \frac{1}{g(x,x)} < \infty. 
$$
This shows that $Q(\Omega^1_M) = 1$, and thus $Q(\Omega^1) = 1$. For every $\omega \in \Omega$, we have 
$$
Z_{x,s}^{y,t}(\omega) > 0, \quad \forall x, y \in \Z^d, \ \forall s, t \in \R, s < t. 
$$
Moreover, as shown above, for every $M \in \N$ and for all $x, y \in \Z^d$, 
$$
\left \langle \sup_{s, t \in (-M,M), s < t} Z_{x,s}^{y,t} \right \rangle  < \infty. 
$$
As a result, we have for $Q$-almost every $\omega \in \Omega$ 
$$
0 < Z_{x,s}^{y,t}(\omega) < \infty \quad \forall x, y \in \Z^d, \ s, t \in \R, s < t. 
$$
By subtracting from $\Omega^1$ a set of measure $0$ and calling the resulting set still $\Omega^1$, we can then assume without loss of generality that for every $\omega \in \Omega^1$ 
\begin{equation}    \label{eq:trivial_E_bounds} 
0 < Z_{x,s}^{y,t}(\omega) < \infty 
\end{equation}
for all $x, y \in \Z^d$ and $s, t \in \R$ such that $s < t$.  

Fix $\omega \in \Omega^1$, $f \in \Lc$, and $s, t \in \R$ such that $s < t$. We need to show that $L^{s,t}_{\omega} f \in \Lc$, i.e., we need to show that there exist $\tilde c > 0$ and $\tilde \epsilon \in (0,1)$, depending on $\omega$, $f$, $s$, and $t$, such that 
$$
e^{-\tilde c \|y\|^{1-\tilde \epsilon}} \leq L^{s,t}_{\omega} f(y) \leq e^{\tilde c \|y\|^{1-\tilde \epsilon}}, \quad \forall y \in \Z^d. 
$$
Since $f \in \Lc$, there exist $c > 0$ and $\epsilon \in (0,1)$ such that $e^{-c \|x\|^{1-\epsilon}} \leq f(x) \leq e^{c \|x\|^{1-\epsilon}}$ for all $x \in \Z^d$. Let $M \in \N$ be so large that $s, t \in (-M,M)$. Since $\omega \in \Omega^1_M$, we can estimate the numerator of $L^{s,t}_{\omega} f(y)$ for $y \in \Z^d$ such that $\|y\| > R_M(\omega)$ as follows: 
\begin{equation}    \label{eq:closed_cocycle_1}  
\sum_{x \in \Z^d} f(x) Z_{x,s}^{y,t}(\omega) \leq \sum_{x \in \Z^d} e^{c \|x\|^{1-\epsilon}} (p_{2M}^{y-x})^{\frac{1}{2}} g(x,y). 
\end{equation}
Estimating $e^{-2M}$ from above by $1$, one obtains  
\begin{align*} 
p_{2M}^{y-x} =& e^{-2M} \sum_{n=0}^{\infty} \frac{(2M)^n}{n!} q_n^{y-x}  < \sum_{n=0}^{\infty} \frac{(2M)^n}{n!} q_n^{y-x} \\
<& \sum_{n=0}^{\infty} \frac{(2M)^n}{n!} q_n^{y-x} + \sum_{n=0}^{\infty} \sum_{m \neq n} \left(\frac{(2M)^n}{n!} q_n^{y-x} \right)^{\frac{1}{2}} \left(\frac{(2M)^m}{m!} q_m^{y-x} \right)^{\frac{1}{2}} 
\end{align*} 
and thus  
$$
(p_{2M}^{y-x})^{\frac{1}{2}} < \sum_{n=0}^{\infty} \left(\frac{(2M)^n}{n!} q_n^{y-x} \right)^{\frac{1}{2}}. 
$$
The transition probability $q_n^{y-x}$ can only be positive if $\|y-x\|_1 \leq n$. If $\|x\| > \|y\| + n$, we have 
$$
\|y-x\|_1 \geq \|y-x\| \geq \|x\| - \|y\| > n, 
$$
so $q_n^{y-x} = 0$. It follows that the expression on the right-hand side of~\eqref{eq:closed_cocycle_1} is bounded from above by 
\begin{equation}   \label{eq:closed_cocycle_2} 
(1+\|y\|)^{d+1} \sum_{n=0}^{\infty} \left(\frac{(2M)^n}{n!} \right)^{\frac{1}{2}} Y_n, 
\end{equation} 
where 
$$
Y_n := \sum_{x: \|x\| \leq \|y\|+n} \left(q_n^{y-x} \right)^{\frac{1}{2}} e^{c \|x\|^{1-\epsilon}} (1+\|x\|)^{d+1}.  
$$
We split the expression from~\eqref{eq:closed_cocycle_2} in two: 
\begin{align}   
& (1+\|y\|)^{d+1} \sum_{0 \leq n \leq \|y\|} \left(\frac{(2M)^n}{n!} \right)^{\frac{1}{2}} Y_n \label{eq:closed_cocycle_3} \\
+& (1+\|y\|)^{d+1} \sum_{n > \|y\|} \left(\frac{(2M)^n}{n!} \right)^{\frac{1}{2}} Y_n.   \label{eq:closed_cocycle_4} 
\end{align}
For $0 \leq n \leq \|y\|$, we have the estimate 
\begin{align*}
Y_n \leq& \sum_{x: \|x\| \leq 2 \|y\|} \left(q_n^{y-x} \right)^{\frac{1}{2}} e^{c \|x\|^{1-\epsilon}} (1+\|x\|)^{d+1} \\
\leq& e^{c 2^{1-\epsilon} \|y\|^{1-\epsilon}} (1+2\|y\|)^{d+1} \sum_{x: \|x\| \leq 2 \|y\|} \left(q_n^{y-x} \right)^{\frac{1}{2}} \\
\lesssim& e^{c 2^{1-\epsilon} \|y\|^{1-\epsilon}} (1+2 \|y\|)^{d+1} \|y\|^d.  
\end{align*}
If we set 
$$
m(a) := (1+a)^{d+1} (1+2a)^{d+1} a^d, 
$$
we can bound the expression in~\eqref{eq:closed_cocycle_3} from above by a constant times  
$$
m(\|y\|) e^{c 2^{1-\epsilon} \|y\|^{1-\epsilon}} \sum_{n=0}^{\infty} \left(\frac{(2M)^n}{n!} \right)^{\frac{1}{2}}. 
$$
Suppose now that $n > \|y\|$. Then  
\begin{align*}
Y_n \leq& \sum_{x: \|x\| \leq 2n} \left(q_n^{y-x} \right)^{\frac{1}{2}} e^{c \|x\|^{1-\epsilon}} (1+\|x\|)^{d+1} \\
\leq& e^{c 2^{1-\epsilon} n^{1-\epsilon}} (1+2n)^{d+1} \sum_{x: \|x\| \leq 2n} \left(q_n^{y-x} \right)^{\frac{1}{2}} \\
\lesssim& e^{c 2^{1-\epsilon} n^{1-\epsilon}} (1+2n)^{d+1} n^d. 
\end{align*}
The expression in~\eqref{eq:closed_cocycle_4} is hence bounded from above by a constant times 
$$
e^{c 2^{1-\epsilon} n^{1-\epsilon}} (1+2n)^{d+1} n^d (1+\|y\|)^{d+1} \left[\sum_{n=0}^{\infty} \left(\frac{(2M)^n}{n!} \right)^{\frac{1}{2}} \right], 
$$
where the series in brackets converges. We have thus shown that there exists a polynomial $p$, with coefficients depending only on $M, \beta$, and $d$, such that 
$$
\sum_{x \in \Z^d} f(x) Z_{x,s}^{y,t}(\omega) \leq p(\|y\|) e^{c 2^{1-\epsilon} \|y\|^{1-\epsilon}}
$$
for all $y \in \Z^d$ with $\|y\| > R_M(\omega)$. In light of~\eqref{eq:trivial_E_bounds}, it is then possible to choose $\tilde c_1 > 0$ so large that 
$$
\sum_{x \in \Z^d} f(x) Z_{x,s}^{y,t}(\omega) \leq \sum_{x \in \Z^d} f(x) Z_{x,s}^{0,t}(\omega) e^{\tilde c_1 \|y\|^{1-\epsilon}}, \quad \forall y \in \Z^d, 
$$
and thus 
$$
L^{s,t}_{\omega} f(y) \leq e^{\tilde c_1 \|y\|^{1-\epsilon}}, \quad \forall y \in \Z^d. 
$$
The proof of the lower bound for $L^{s,t}_{\omega} f$ is simpler: Since $\omega \in \Omega^1_M$, we have for every $y \in \Z^d$ with $\|y\| > R_M(\omega)$ 
\begin{equation*}  
\sum_{x \in \Z^d} f(x) Z_{x,s}^{y,t}(\omega) \geq \sum_{x \in \Z^d} e^{-c \|x\|^{1-\epsilon}} Z_{x,s}^{y,t}(\omega) \geq e^{-c \|y\|^{1-\epsilon}} Z_{y,s}^{y,t}(\omega)
\geq \frac{e^{-c \|y\|^{1-\epsilon}}}{g(y,y)}. 
\end{equation*}
Then, again by virtue of~\eqref{eq:trivial_E_bounds}, we can choose $\tilde c_2 > 0$ so large that 
$$
\sum_{x \in \Z^d} f(x) Z_{x,s}^{y,t}(\omega) \geq \sum_{x \in \Z^d} f(x) Z_{x,s}^{0,t}(\omega) e^{-\tilde c_2 \|y\|^{1-\epsilon}}, \quad \forall y \in \Z^d, 
$$
and hence we have for $\tilde c:= \max\{\tilde c_1, \tilde c_2\}$ the estimate 
$$
e^{-\tilde c \|y\|^{1-\epsilon}} \leq L^{s,t}_{\omega} f(y) \leq e^{\tilde c \|y\|^{1-\epsilon}}, \quad \forall y \in \Z^d. 
$$

\bigskip 

\subsection{Proof of Lemma~\ref{lm:Omega_subset}} 

Let us first show that there exists a set $\Omega^Z \in \Fc$ such that $Q(\Omega^Z) = 1$ and the function 
$$
\mathcal{Y}(\omega): \Z^d \to (0, \infty), \ y \mapsto Z_{-\infty}^{y,0}(\omega)/Z_{-\infty}^{0,0}(\omega)
$$
is an element of $\Lc$ for every $\omega \in \Omega^Z$. Recall from Proposition~\ref{prop:global_stationary_sol} that for $\beta$ sufficiently small there is $\Omega^{\textup{sol}} \in \Fc$ such that $Q(\Omega^{\textup{sol}}) = 1$ and such that 
$$
\begin{aligned}
	\Z^d \times \R \times \Omega^{\textup{sol}} 
		&\to (0, \infty)
\\	(y,t,\omega)
		&\mapsto Z_{-\infty}^{y,t}(\omega)
\end{aligned}
$$
is a global stationary solution to~\eqref{eq:sSHE}. Let $\Omega^Z$ be the set given by 
\begin{equation}
\Omega^Z  
	:=
		\left\{\omega \in \Omega^{\textup{sol}} ~:~
			\substack{ \displaystyle 
				\ \exists R(\omega) > 0 \text{ s.t } \forall y\in \Z^d \text{ with } \|y\| > R(\omega):
			\\	\displaystyle
				(1+\|y\|)^{-(d+1)} \leq Z_{-\infty}^{y,0}(\omega) \leq (1+\|y\|)^{d+1}
			}
			\right\}.
\end{equation}
For every $\omega \in \Omega^Z$, the function $\mathcal{Y}(\omega)$ has polynomial growth and is in particular an element of $\Lc$. 

For $y \in \Z^d$, define the sets 
\begin{align*}
E_1(y) :=& \{\omega \in \Omega^{\textup{sol}}: \ Z_{-\infty}^{y,0}(\omega) > (1+\|y\|)^{d+1}\}, \\
E_2(y) :=& \{\omega \in \Omega^{\textup{sol}}: \ Z_{-\infty}^{y,0}(\omega) < (1+\|y\|)^{-(d+1)}\}, 
\end{align*}
and let $E(y) := E_1(y) \cup E_2(y)$. 
By the Borel--Cantelli lemma, in order to show that $Q(\Omega^Z) = 1$, it is enough to prove that $\sum_{y \in \Z^d} Q(E(y)) < \infty$. By Markov's inequality,  
$$
Q(E_1(y)) \leq \frac{\langle Z_{-\infty}^{y,0} \rangle}{(1+\|y\|)^{d+1}} = (1+\|y\|)^{-(d+1)}. 
$$
Similarly, 
$$
Q(E_2(y)) \leq \frac{\langle (Z_{-\infty}^{y,0})^{-1} \rangle}{(1+\|y\|)^{d+1}} = \frac{C}{(1+\|y\|)^{d+1}} 
$$
for some $C > 0$ (see Theorem~\ref{thm:continuous_Talagrand} and the subsequent remark on negative moments for the limiting partition function).  Hence, 
$$
\sum_{y \in \Z^d} Q(E(y)) \leq \sum_{y \in \Z^d} \frac{C+1}{(1+\|y\|)^{d+1}} < \infty. 
$$
Now we use $\Omega^Z$ to construct a new set $\widetilde{\Omega}$ that has the desired properties. By Lemma~\ref{lm:L_invariance}, there exists a set $\Omega^L \in \Fc$ such that $Q(\Omega^L) = 1$ and $L_{\omega}^{s,t} f \in \Lc$ for every $f \in \Lc$, $s, t \in \R$ with $s \leq t$, and $\omega \in \Omega^L$. Set 
$$
\Omega_1 := \bigcap_{r \in \Z} \theta_r(\Omega^Z \cap \Omega^L). 
$$
We claim that for every $\omega \in \Omega_1$ and for every $t \in \R$, the function $y \mapsto \widetilde{Z}_{-\infty}^{y,t}(\omega) := Z_{-\infty}^{y,t}(\omega)/Z_{-\infty}^{0,t}(\omega)$ is an element of $\Lc$. First consider the case $t > 0$. For $\omega \in \Omega_1$ and $y \in \Z^d$, one has 
$$ 
\widetilde{Z}_{-\infty}^{y,t}(\omega) = \frac{\sum_{x \in \Z^d} Z_{-\infty}^{x,0}(\omega) Z_{x,0}^{y,t}(\omega)}{\sum_{x \in \Z^d} Z_{-\infty}^{x,0}(\omega) Z_{x,0}^{0,t}(\omega)} = L_{\omega}^{0,t} \mathcal{Y}(\omega).  
$$ 
Since $\mathcal{Y}(\omega) \in \Lc$ and since $\omega \in \Omega^L$, it follows that $y \mapsto \widetilde{Z}_{-\infty}^{y,t}(\omega)$ is an element of $\Lc$. Now assume that $t \leq 0$ and let $r \in \Z$ be so large that $t+r > 0$. For $\omega \in \Omega_1$ and $y \in \Z^d$, one has 
$$
\widetilde{Z}_{-\infty}^{y,t}(\omega) =  \widetilde{Z}_{-\infty}^{y,t+r}(\theta_{-r} \omega)
$$
and $y \mapsto \widetilde{Z}_{-\infty}^{y,t+r}(\theta_{-r} \omega)$ is an element of $\Lc$. This proves the claim. 

Consider the set 
\begin{align*}
\widetilde{\Omega} := \{\omega \in \Omega^{\textup{sol}}: \ & L_{\omega}^{s,t} f \in \Lc \ \forall s, t \in \R \ \text{with} \ s \leq t \ \text{and} \ \forall f \in \Lc; \\
& y \mapsto \widetilde{Z}_{-\infty}^{y,t}(\omega) \in \Lc \ \forall t \in \R\}. 
\end{align*} 
Then $\Omega_1 \subset \widetilde{\Omega}$. As $Q(\Omega_1) = 1$ and as $(\Omega, \Fc, Q)$ is complete, $\widetilde{\Omega}$ is measurable and has measure $1$. The fact that $\widetilde{\Omega}$ is invariant under $\theta_s$ for every $s \in \R$ is seen immediately and the remaining properties are direct consequences of how $\widetilde{\Omega}$ is defined.

\subsection{Proof of Lemma~\ref{lm:cocycle}} 

Proving measurability of $\varphi$ is not difficult and will be omitted. Fix $\omega \in \widetilde{\Omega}$. Then, for $f \in \Lc$ and $y \in \Z^d$,  
$$
\varphi^0_{\omega} f(y) = L^{0,0}_{\omega} f(y) = \frac{f(y)}{f(0)} = f(y). 
$$
It remains to show that 
$$
\varphi^{s+t}_{\omega} = \varphi^t_{\theta_s \omega} \circ \varphi^s_{\omega}, \quad \forall s, t \geq 0. 
$$ 
If $s=0$ or $t=0$, this follows immediately. If $s, t > 0$, one has for $f \in \Lc$ and $y \in \Z^d$ 
\begin{align*}
\left(\varphi^t_{\theta_s \omega} \circ \varphi^s_{\omega} \right) f(y) =& \frac{\sum_{x \in \Z^d} \varphi^s_{\omega} f(x) Z_{x,0}^{y,t}(\theta_s \omega)}{\sum_{x \in \Z^d} \varphi^s_{\omega} f(x) Z_{x,0}^{0,t}(\theta_s \omega)} \\
=& \frac{\sum_{x \in \Z^d} \sum_{z \in \Z^d} f(z) Z_{z,0}^{x,s}(\omega) Z_{x,0}^{y,t}(\theta_s \omega)}{\sum_{x \in \Z^d} \sum_{z \in \Z^d} f(z) Z_{z,0}^{x,s}(\omega) Z_{x,0}^{0,t}(\theta_s \omega)} \\
=& \frac{\sum_{z \in \Z^d} f(z) \sum_{x \in \Z^d} Z_{z,0}^{x,s}(\omega) Z_{x,s}^{y,t+s}(\omega)}{\sum_{z \in \Z^d} f(z) \sum_{x \in \Z^d} Z_{z,0}^{x,s}(\omega) Z_{x,s}^{0,t+s}(\omega)} \\
=& \frac{\sum_{z \in \Z^d} f(z) Z_{z,0}^{y,t+s}(\omega)}{\sum_{z \in \Z^d} f(z) Z_{z,0}^{0,t+s}(\omega)} = \varphi^{t+s}_{\omega} f(y).  
\end{align*}

\bibliographystyle{alpha}

\bibliography{references_v2}

\end{document}